\documentclass{amsart}
\usepackage{graphicx,lscape,mathtools}
\usepackage{amsmath, amsfonts, amssymb, amsthm}
\usepackage{MnSymbol}
\usepackage{ifthen}
\usepackage{caption}
\usepackage{subcaption}
\usepackage[justification=centering]{caption}
\usepackage[colorlinks=true]{hyperref}
\usepackage[export]{adjustbox}
\usepackage{xcolor}
\usepackage{soul}

\newtheorem{theorem}{Theorem}[section]
\newtheorem{corollary}[theorem]{Corollary}
\newtheorem{lemma}[theorem]{Lemma}
\newtheorem{proposition}[theorem]{Proposition}

\theoremstyle{definition}
\newtheorem{definition}[theorem]{Definition}

\newtheorem{example}[theorem]{Example}
\newtheorem{notation}[theorem]{Notation}
\newtheorem{assumptions}[theorem]{Assumptions}

\newcommand{\CC}{\mathbb{C}}

\newcommand{\ZZ}{\mathbb{Z}}
\newcommand{\bA}{\mathbb{A}}

\newcommand{\RR}{\mathbb{R}}

\newcommand{\cM}{\mathcal{M}}

\newcommand{\mandel}{\mathcal{M}}
\newcommand{\Jplus}{J_+}
\newcommand{\Jstar}{J_*}

\newcommand{\Jminus}{J_-}
\newcommand{\Jnot}{J_0}
\newcommand{\Kplus}{K_+}

\newcommand{\Kminus}{K_-}
\newcommand{\Knot}{K_0}
\newcommand{\FC}{\mathcal{U}}

\newcommand{\Rnab}{F_{n,a,b}}
\newcommand{\mapR}{F}
\newcommand{\mapRnab}{F}

\newcommand{\Rna}{r_{n,a}}

\newcommand{\spire}{\mathcal{S}}

\newcommand{\FClabel}[4]{\begin{pmatrix}
    #2 \curvearrowleft #1\\
    #3 \rcurvearrowright #4 \end{pmatrix}}

\DeclareMathOperator\Arg{Arg}

\begin{document}

\title[Paper Fortune Tellers: Sidecars and Zippers]{Paper Fortune Tellers in Julia sets of Generalized McMullen maps II: \\Sidecars and Zippers}

\author[S.~Boyd]{Suzanne Boyd}
\address{Department of Mathematical Sciences\\
University of Wisconsin Milwaukee\\
PO Box 413\\
Milwaukee, WI 53201, 
USA}
\email{sboyd@uwm.edu, ORCID: 0000-0002-9480-4848}

\author[K.~Brouwer]{Kelsey Brouwer}
\address{Mathematics Department\\
Aquinas College\\
1700 Fulton St. E. \\
Grand Rapids, MI 49506\\
USA}
\email{brouwerk@aquinas.edu}

\date{\today}

\begin{abstract}
 We study the family of complex rational functions known as Generalized McMullen maps,  
    \noindent \mbox{$\Rnab(z) = z^n + \dfrac{a}{z^n}+b$,} 
    for $a\neq 0$ and $n \geq 3$ fixed. In \cite{BoydBrouwer1}, we provided a combinatorial model for a large class of maps whose Julia sets contain both infinitely many homeomorphic copies of quadratic Julia sets conjugate to the ``basilica'', and infinitely many subsets  homeomorphic to a set which is obtained by starting with the basilica, then changing a finite number of pairs of external ray landing point identifications, following an algorithm we described.  

In this article, we generalize beyond the basilica, and provide a catalog of additional types of hyperbolic Julia sets of Generalized McMullen maps, where the ``baby'' Julia set can be any rabbit, aeroplane, or Kokopelli quadratic Julia set; that is, where the $c$-value can be taken from any bulb attached to the main cardioid of the Mandelbrot set, or from the main cardioid of any principal baby Mandelbrot set (no renormalizations). 
\end{abstract}

\maketitle

\markboth{\textsc{S. Boyd and K. Brouwer}}
  {\textit{Sidecars and Zippers in Generalized McMullen Julia Sets}}

\footnotetext[1]{2020 MSC-class: 37F10 (Primary) 37F12, 37F20 (Secondary).}

\section{Introduction}

The goal of this article is to describe and catalog some interesting dynamical behavior in Julia sets of the family of complex, rational functions: $$\Rnab(z) = z^n + \dfrac{a}{z^n}+b~,~n \in \mathbb{N},~a \in \mathbb{C}^*,~b \in \mathbb{C},$$ where $\mathbb{C}^* = \mathbb{C}\setminus \lbrace 0\rbrace$. We consider $n$ a fixed integer with $n \geq 3$, thus have two complex parameters $a\neq 0$ and $b$. This family has $2n$ critical points, but only two critical values, we refer to as $v_+$ and $v_-$.  Following Xiao, Qiu, and Yin (\cite{xiaoqiu}), we call this family ``Generalized McMullen Maps'', as McMullen introduced the study of the subfamily $z\mapsto z^n + \dfrac{a}{z^n}$ (\cite{mcmullen_example}).

McMullen maps and Generalized McMullen maps, sometimes considering the alternate generalization $z\mapsto z^n +\dfrac{a}{z^d}$, have also been studied previously by Devaney and colleagues (\cite{DevaneyHalos, DevaneySurvey2013, DevaneyRussell, DevaneyKozma, devgar, DevaneyHalos, jangso} and see \cite{Stoertz-error}),  as well as Boyd and colleagues (\cite{BoydSchulz, BoydMitchell, BoydHoeppner1, BoydBrouwer1, BoydBrouwerHoeppner}).  Stoertz and colleagues have begun the study of ``Maximally Generlized McMullen Maps'' $z\mapsto z^n + b +\dfrac{a}{z^d}$ (\cite{Murali-PME}).

\begin{figure}\centering
\includegraphics[width=.5\textwidth]{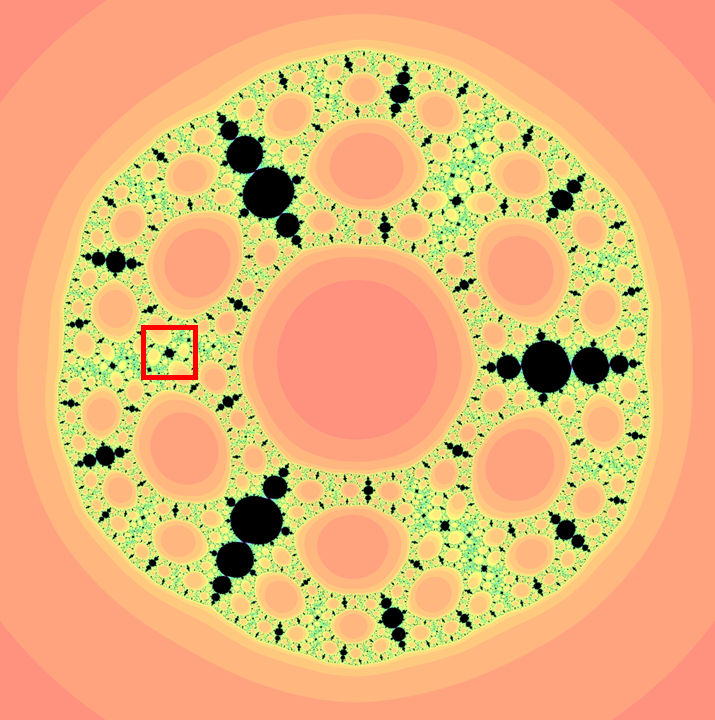}
    \caption[A rational Julia set with standard and altered Julia sets]{\label{fig:wholeJexample}
    The Julia set of $\Rnab$ for $n=3$, $a=0.05855-0.01282i$, $b=0.02+0.03i$. The baby quadratic Julia set is on the positive real axis. Its $n-1=2$ rotationally symmetric preimages are apparent, as are several smaller deeper-level preimages. An ``altered'' preimage is several preimages deep, so it is too small to see in detail without zooming in to the area outlined by the red square.}
\end{figure}

These maps, even for $b=0$, exhibit some behavior that is distinct from polynomial dynamics; for example, there are Julia sets which are Cantor sets of simple closed curves (\cite{mcmullen_example, devlook}. However, many Julia sets in this family contain homeomorphic copies of quadratic polynomial Julia sets, ``baby Julia sets'' as named by Douady and Hubbard (\cite{DouadyHubbard}), strewn throughout a ``necklace'' structure - see Figure~\ref{fig:wholeJexample}. The form of the map guarantees an $n$-fold rotational symmetry in the Julia set.

In prior work by the authors of the present article \cite{BoydBrouwer1}, we reveal and analyze the dynamical behavior in the case that one of the (two) critical orbits of the map is responsible for a baby quadratic Julia set which is homeomorphic to a basilica (the Julia set for $z\mapsto z^2-1)$, and the other critical value (strictly eventually) lands in that same baby basilica (hence, in the language of Milnor, a ``capture'' situation). Our first step in doing so is defining a set of external angles (a map from $[0,1)$) onto the baby Julia set (in the way preserving the expected dynamics), and then pulling back to define external angles on every one of its preimages, in the rational Julia set, in a way that respects the rational dynamics. (See Theorem~\ref{thm: angle assignments}.)
We show the consequence of having the second critical value in the (eventual) preimage of a baby basilica, located elsewhere in the rational map's Julia set, is that if you pull back to one preimage before the other critical value, you see a shape which in many ways resembles the baby basilica, \textit{except} it has a finite number of external angle identification changes--combining some Fatou components while splitting others. We call this an ``altered'' basilica. 

\begin{figure}
    \centering
    \begin{minipage}{0.455\textwidth}
        \includegraphics[width=\textwidth]{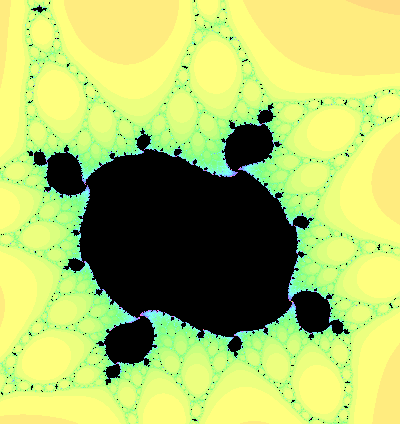}
    \end{minipage}
    \begin{minipage}{0.445\textwidth}
        \includegraphics[width=\textwidth]{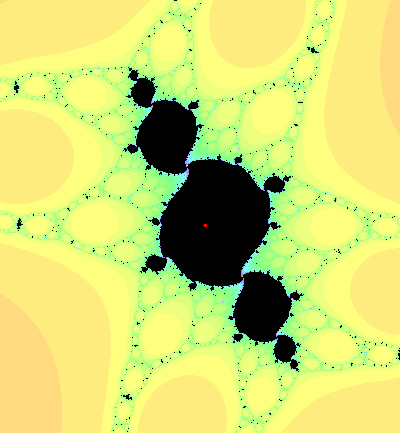}
    \end{minipage}
    \caption[An ``altered'' and a standard Julia set]{\label{fig:introexample} Portions of the Julia set of $\Rnab$ for $n=5$, $a=0.1317-0.0073i$, $b=0.03+0.02i$. The right image is homeomorphic to a quadratic Julia set, as it is a preimage copy under $\Rnab$ of a baby quadratic Julia set in $J(\Rnab)$ associated with the critical value $v_+$. The left is a preimage of the right under $\Rnab$, and is what we refer to as ``altered''. The red dot in the center of the right image marks the location of the other critical value $v_-$; note this is \textit{not} in the standard location of a critical value in a basilica Julia set, which would be in one of the largest Fatou components adjacent to the central component. This alternate placement causes the altered preimage shape. }
\end{figure}

\begin{figure}
  \includegraphics[width=.375\textwidth]{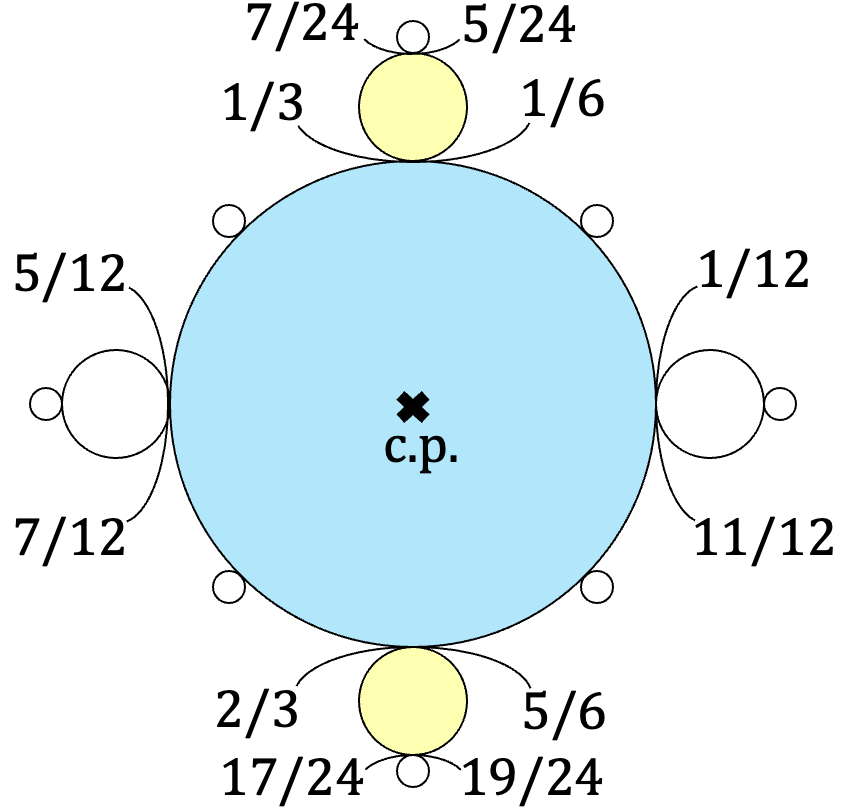}
    \includegraphics[width=.475\textwidth]{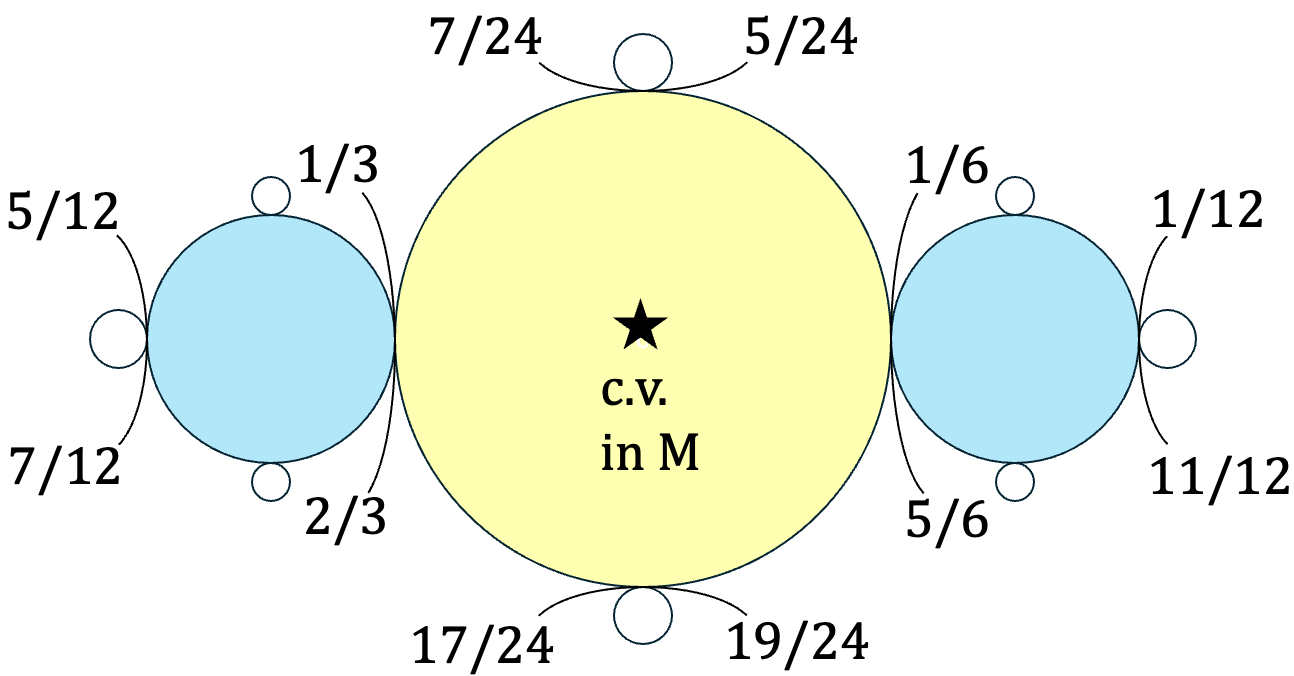}
\caption{\label{fig:vminus-in-M-bas}
On the right is a schematic of the right of Figure~\ref{fig:introexample} (a preimage of a baby basilica), with the other critical value in the central component labelled ``M''. On the left is a schematic of the left of Figure~\ref{fig:introexample}, the preimage of the right.} 
          \end{figure}

In \cite{BoydBrouwer1}, we provide a combinatorial model description of this process via taking some sets of rays that are identified (i.e., land at the same Julia set point), breaking the identifications, then gluing them into different pairings. This can be visualized somewhat as holding an origami paper ``fortune teller'' open in one direction, then closing it and opening it in the other.   
We show precisely which ray identifications change depending upon the location of the other critical value in the preimage of the baby basilica.  See Figure~\ref{fig:introexample} for a baby basilica and its ``altered'' preimage. 
See Figure~\ref{fig:vminus-in-M-bas} for a schematic of the situation in Figure~\ref{fig:introexample}, showing the key angle identification changes on the preimage of the baby $J$. These changes primarily occur because the preimage of a critical value component can only be one Fatou component which must map 2:1 onto its image. 

In this article, we expand that study significantly by exploring further in the parameter space of the Generalized McMullen maps, now allowing the baby Julia set to be any hyperbolic, connected, quadratic Julia set with parameters taken from a rabbit bulb (primary decoration of the main cardioid) or from an aeroplane or Kokopelli cardioid (main cardioid of any ``top-level'' baby Mandelbrot set).
We include a selection of representative examples to give the reader a good understanding of the effects of both the type of baby Julia set and the location of the other critical value within the preimage copy of the baby Julia set. As we will see, $\nu$-rabbits cause ``sidecars'', and aeroplanes require visualizing a ``zipper'' change in angle identifications.
We note there are additional examples in the PhD thesis of the second author (\cite{Brouwer-Thesis}).

Our main theorems (Theorems~\ref{thm: rabbit type N} and~\ref{thm: zip type N}) require detailed notation and definitions, and are best understood after working through the examples we present, so we do not state them in this introduction. 

This article is organized as follows.  Section~\ref{sec:prelim} reviews needed background and preliminaries, including from \cite{BoydBrouwer1},
and 
briefly describes the somewhat degenerate case in which the baby Julia set is a topological disk (thus there are no ray identifications to change, but there is a still an alteration in a sense). 
Section~\ref{sec:sidecars} contains the description of the altered baby Julia sets in the case of any $\nu$-rabbit, coming from a bulb attached to the main cardioid of the Mandelbrot set, $\mandel$, including making a connection to the work in \cite{BoydBrouwer1} by thinking of the basilica as a $2$-rabbit (with a ``sidecar'' as we explain).
Section~\ref{sec:zippers} presents the aeroplane and Kokopelli baby Julia set alterations, which are quadratic maps which lie in a main cardioid of a principal baby Mandelbrot set within the Mandelbrot set (this requires introducing ``zippers''). 
For the reader who wishes to jump to see an exotic altered baby Julia set, we suggest 
Figure~\ref{fig: koko N=3}.

\subsection*{Acknowledgements} We appreciate Danny Stoll for his computer program ``Dynamo'' which we used to identify locations of external rays in quadratic Julia sets. We also appreciate Brian Boyd for use of the program ``Dynamics Explorer'' to generate parameter and dynamical space images in this article.

\section{Preliminaries and Background}
\label{sec:prelim}

As ``baby'' Julia sets, the quadratic polynomials $P_c(z) = z^2 + c,~c \in \mathbb{C}$ serve in a way as building blocks for our Generalized McMullen Julia sets. Recall the \textbf{Fatou set} of a polynomial or rational is the set of $z$-values in the domain with stable behavior, that is, where the iterates form a normal family in the sense of Montel,  and the \textbf{Julia set} $J$ is the complement to the Fatou set. For any map where $\infty$ is a super-attracting fixed point (true for polynomials and our Generalized McMullen maps), the \textbf{filled Julia set} $K$ is the union of the Julia set and the bounded Fatou components, so that $K$ is the set of all points with bounded orbits and $J=\partial K$. For the family $P_c$, the \textbf{Mandelbrot set} $\mandel$ can be defined as the set of parameters $c$ with connected filled Julia sets $K_c$.

\subsection{Our maps of interest}

Douady and Hubbard (\cite{DouadyHubbard}) showed that quadratic Julia sets can result from other iterative processes, which explains the existence of homeomorphic copies of polynomial Julia sets inside of other Julia sets. 

\begin{definition}(\cite{DouadyHubbard})
\label{defn:polynlike}
A map $f:U' \to f(U')=U$ is \textbf{polynomial-like} if 
\begin{itemize}
    \item $U'$ and $U$ are bounded, open, simply connected subsets of $\CC$,
    \item $U'$ is relatively compact in $U$, and
    \item $f$ is analytic and proper.
\end{itemize}
Further, $f$ is polynomial-like of \textbf{degree two} if $f$ is a 2-to-1 map except at finitely many points, and $U'$ contains a unique critical point of $f$. The \textbf{filled Julia set of a polynomial-like map} is the set of points whose orbits remain in $U'$. 
\end{definition}

\begin{theorem} \label{thm:polynlike}
    \cite{DouadyHubbard} A polynomial-like map of degree two is topologically conjugate on its filled Julia set to a quadratic polynomial on that polynomial's filled Julia set.
\end{theorem}

For this reason, the filled Julia set of a polynomial-like map is often referred to as a \textbf{baby Julia set}. 
Our results require that $\Rnab$ contains a baby Julia set.
More precisely, in this article, we focus on the maps which satisfy the following set of Assumptions.

\begin{assumptions} \label{ass:general}
Consider $\Rnab$ with integer $n\geq 3$, $a \in \CC^{*}$, $b \in \mathbb{C}$ and for which:
\begin{enumerate}
    \item[(A1)] $\mapR$ is hyperbolic (i.e., the parameters are in a hyperbolic component of parameter space),
    \item[(A2)] Both critical values have bounded orbits (i.e., $v_+, v_-$ of $\Rnab$ lie in the filled Julia set $K_{n,a,b} = K(\Rnab)$), but are not in the same Fatou component,
    \item[(A3)] $\mapR$ is polynomial-like of degree 2 on a region $U'$ containing one critical value (wolog $v_+$), and conjugate on the filled Julia set of $\mapR|_{U'}$, which we'll call $\Kplus$, to a quadratic polynomial $P_c$,
    \item[(A4)] $P_c$ is connected and hyperbolic; thus $v_+$ lies in the immediate basin of a finite attracting cycle, and finally,
    \item[(A5)] The other critical value $v_-$ lies in the interior of a preimage copy of the baby quadratic Julia set $K_+$,  and therefore lies in the same attracting basin as $v_+$ but not in the immediate basin. We call that preimage $K_-$. 
\end{enumerate}
\end{assumptions}

In \cite{BoydBrouwer1}, Assumption (A4) above was restricted to $P_c$ in the ``basilica'' bulb of the Mandlebrot set, in which the critical orbit is periodic of period $2$. As stated earlier, in our results in this article, we explore some cases for which $P_c$ is different types of baby Julia sets.

\subsection{External Angle Assignments}

We first recall the definition of external rays and angles as described in \cite{milnor}. 
For a simply connected filled Julia set $K(P)$, for a polynomial map $P$ of degree $n\geq 2$,
the Böttcher coordinate $\phi: \CC \smallsetminus K \to \CC \smallsetminus \mathbb{D}$ is an isomorphism conjugating $P$ to $z \mapsto z^n$. The sets $\{z : \Arg(\phi(z)) \text{ is constant}\}$ generate \textbf{external rays} for $K$. Alternatively, we denote by $\rho_t$ the ray $\{\phi^{-1}(re^{2\pi it}): r>1\}$,  defined on $\CC\setminus K$. If $\gamma(t) = \displaystyle \lim_{r\searrow1} \phi^{-1}(re^{2\pi it})$ exists, the ray \textbf{lands} and we can associate $z=\gamma(t)$ with the \textbf{external angle} $t$.

All rays land whenever $J$ or $K$ is locally connected, which is true for our maps of interest. Further, if the angle $t\in \mathbb{R} \slash \mathbb{Z}$ is rational, the landing point $\gamma(t)$ is either a periodic or pre-periodic point. We  identify $S^1$ with the unit interval $[0,1)$.
See Figure~\ref{fig:labels} for a quadratic polynomial Julia set (the ``basilica'') with some labeled rays. 
Note multiple rays often land at the same point.

Similarly, \cite{milnor} defines \textbf{parameter rays} for the outside of the Mandelbrot set, $\cM$, where rays of rational angles land on a point in $\partial \cM$. 

We use the following result of \cite{BoydBrouwer1} to assign sets of external angles to the subset of the Julia set consisting of $\Jplus=\partial \Kplus$ and its tree of preimages. Though most of \cite{BoydBrouwer1} focused on basilicas, we kept the following result and its proof applicable to any Generalized McMullen map satisfying Assumptions~\ref{ass:general}.

\begin{theorem}[\cite{BoydBrouwer1}] \label{thm: angle assignments}
    Let $F=\Rnab$ be a generalized McMullen map that is polynomial-like on a region $U'$ containing the critical value $v_+$, and that is conjugate on $K_+$, the filled Julia set of $F|_{U'}$, to a quadratic polynomial $P_c$ on its filled Julia set $K(P_c)$. Assume that $K(P_c)$ is connected and locally connected. Let $\Jstar = \cup_{m=0}^{\infty} F^{-m}(J_+)$.

    Then there exists a surjective relation $\Gamma:S^1\to \Jstar$ which assigns an angle in $[0,1)$ to each point in $\Jstar$ so that the angle assignments respect the dynamics to and from that point. In particular, if $J_{m,j}\in F^{-m}(\Jplus)$ is a preimage copy of $\Kplus$, then $\Gamma$ restricted to co-domain $J_{m,j}$ is a surjective function $\Gamma|^{J_{m,j}}:S^1\to J_{m,j}$, and
    \begin{enumerate}
        \item If $m=0$, $F(\Gamma|^{J_+}(t))=\Gamma|^{J_+}(2t)$.
        \item If $m\geq1$ and $K_{m,j}\neq K_+$ is a component of $F^{-m}(K_+)$ which contains a critical point of $F$, then $F(\Gamma|^{J_{m,j}}(t))=\Gamma|^{F(J_{m,j})}(2t)$.
        \item If $m>1$ and $K_{m,j}$ is a component of $F^{-m}(K_+)$ which does not contain a critical point of $F$, then $F(\Gamma|^{J_{m,j}}(t))=\Gamma|^{F(J_{m,j})}(t)$.
    \end{enumerate}
\end{theorem}

The proof given in \cite{BoydBrouwer1} works inductively by first taking the identification of $S^1$ with the interval $[0,1)$ and pulling those angles back to $J_+$ using the map $\phi$ that conjugates the dynamics of $F$ on $K_+$ to that of $P_c$ on $K(P_c)$. Then those angle assignments can be pulled back through successive preimages of $J_+$ (whose union we called $J_*$) to assign angles on each eventual preimage copy of $J_+$ in a way that respects the dynamics.
Since $F$ is polynomial-like on $K_+$, $K_+$ must contain a critical value, so each of its direct preimages must contain a critical point. The local dynamics of $F$ are degree 2, so we see that the direct preimages of $J_+$ all map $2:1$ onto $J_+$, so the relation $\Gamma$ that assigns angles onto each preimage is conjugate to angle doubling on this domain. At each sequential set of preimages, angles are assigned based on whether the new preimages contain a critical point or not. If not, $\Gamma$ is conjugate to the identity map and angles are assigned as ``clones'' of their images. If a critical point is present, $\Gamma$ is again conjugate to angle doubling. 

While this result from \cite{BoydBrouwer1} will be re-used as written, some other definitions, notations, and results from \cite{BoydBrouwer1} will be slightly altered in this work to generalize better to various baby Julia sets. 
We begin with listing some essential notation.

\begin{notation} \label{not:FC_- etc}
Consistent with the notation from \cite{BoydBrouwer1}, suppose that $v_- \in K_{\ell, j}$, where $K_{\ell,j}$ is a distinct component of $\mapR^{-\ell}(\Kplus)$ for some $\ell \geq 1$ (so $j$ is just an index to distinguish such components). For ease of notation, we use:
\begin{itemize}
\item 
$\Kminus = K_{\ell, j}$ and $\Jminus = \partial \Kminus$, 
so  $v_- \in \Kminus$ parallels $v_+ \in \Kplus$;
\item 
$\FC_-$ is the Fatou component  of $\Kminus$ in which $v_-$ lies;
\item $\FC_c$ is the Fatou component of $K_c$ which corresponds to $\FC_-$;  
\item  $\Jnot=J_{\ell+1, j'}$  refers to any preimage component in $\mapR^{-1}(\Jminus)$, so that $\mapR(J_0)= J_-$; and
\item given a $\Jnot=J_{\ell+1,j'}$, let 
$\gamma_- = \gamma_{\ell,j}:S^1 \to \Jminus$ and 
$\gamma_0 = \gamma_{\ell+1,j'}:S^1 \to \Jnot$ be the shorthand notation for the $\Gamma|^{\Jminus}$ and $\Gamma|^{\Jnot}$ constructed in Theorem~\ref{thm: angle assignments}.
\end{itemize}
\end{notation}

Under Assumptions~\ref{ass:general}, $\Kminus$ is always  a ``clone'' of the baby Julia set $\Kplus$, as in it has the same identified angles and relative size and position of Fatou components. The key for this study is that nothing is forcing $v_-$ to lie in the same relative position in $\Kminus$ as $v_+$ does in $K_+$ (which is the same as $c$ does in $P_c$). A priori, $v_-$ could be anywhere in $\Kminus$, as the two critical orbits are independent. 
This is important, as the position of $v_-$ within $\Kminus$ determines the changes to angle identifications and Fatou components in $\Jnot$ and its tree of preimages as compared to $\Jminus$ and its clones.  

\textbf{Wolog.} 
Note that in \cite{BoydBrouwer1}, we always assume that the critical value $v_-$ lies somewhere on the ``left'' side of $\Kminus$ or that it lies in some smaller decoration attached to the central component of the baby basilica. 
We can make this assumption because 
at the step of the angle assignment construction in which $\Kminus$ receives its angles,
there are always two points of $\Kminus$ to which map to the point with angle zero on $\Kplus$, so we may choose which of these points to assign the angles 0 and $\frac12$ so that $v_-$ lies as close as possible to its expected component location rather than in a symmetric component across $\Kminus$.
This simplifies the analysis into fewer cases. In general, we made this choice, then rotated the image to have the point with angle 0 appear on the right side of the image for the sake of consistency.

\subsection{Main Cardioid Results}
\label{sec:main cardioid}

In this subsection, we provide a result about the case in which in Assumption (A4), $c$ is chosen from the main cardioid so that $K(P_c)$ takes the form of a topological disk. There, no angle identifications exist to complicate matters. 
We note that in a previous work, \cite{BoydBrouwerHoeppner}, the authors of this article along with Hoeppner studied the subfamily $r_{n,a}(z) = z^n + \frac{a}{z^n}+ (a^{1/2n}-2\sqrt{a})$, which uses $b:=a^{1/2n}-2\sqrt{a}$ to force the critical value $v_+$ to be a fixed point. Under this restriction, the baby Julia set in $J(r_{n,a})$ can only be a topological disk, on which $r_{n,a}$ is conjugate to the map $z \mapsto z^2$ on the closed unit disk. There we also applied Theorem~\ref{thm: angle assignments} to obtain a result similar to this one, to apply to baby Julia sets on which $R_{n,a,b}$ is conjugate more generally to $P_c$ on its Julia set, for $c$ chosen from anywhere in the main cardioid of the Mandelbrot set. We do not provide any additional proof as the following is a straightforward application of Theorem~\ref{thm: angle assignments}. 

\begin{proposition} \label{prop: preimage classes}
 Let $\mapR=\Rnab$ be a generalized McMullen map satisfying Assumptions~\ref{ass:general}; moreover, for $P_c$, suppose $c$ is chosen from the main cardioid of the Mandelbrot set $\mandel$.
    
    Then $J(F)$ also contains the infinite tree of preimages $J_* = \cup_{m=0}^\infty\  F^{-m}(\Jplus)$, where each preimage component satisfies one of the following:
    \begin{enumerate}
        \item For each component $J_{1,j}$ of $F^{-1}(\Jplus)$ where $j=0,\dots,n-1$, which includes $\Jplus$ itself, we have\\
        $F\Big(\Gamma|^{J_{1,j}}\big([0,\frac12)\big)\Big) = F\Big(\Gamma|^{J_{1,j}}\big([\frac12,1)\big)\Big) = \Gamma|^{\Jplus}\big([0,1)\big)$.
        
        \item For components $J_{m,j}$ of $F^{-m}(\Jplus)$ such that none of $F^{k}(J_{m,j})$ contain a critical point for $k=0,\dots,m-2$, then we have \\
        $F^m\Big(\Gamma|^{J_{m,j}}\big([0,\frac12)\big)\Big) = F^m\Big(\Gamma|^{J_{m,j}}\big([\frac12,1)\big)\Big) = \Gamma|^{\Jplus}\big([0,1)\big)$.
        
        \item For components $J_{m,j}$ of $F^{-m}(\Jplus)$ such that $F^\ell(J_{m,j})$ contains a critical point for $\ell \in 0,\dots,m-2$, then we have \\
        $F^m\Big(\Gamma|^{J_{m.j}}\big([x,x+\frac14)\big)\Big)=\Gamma|^{\Jplus}\big([0,1)\big)$ for each of $x=0, \frac14, \frac12, \frac34.$
    \end{enumerate}
\end{proposition}

\begin{figure}
  \centering
\includegraphics[width=.5\textwidth,keepaspectratio]{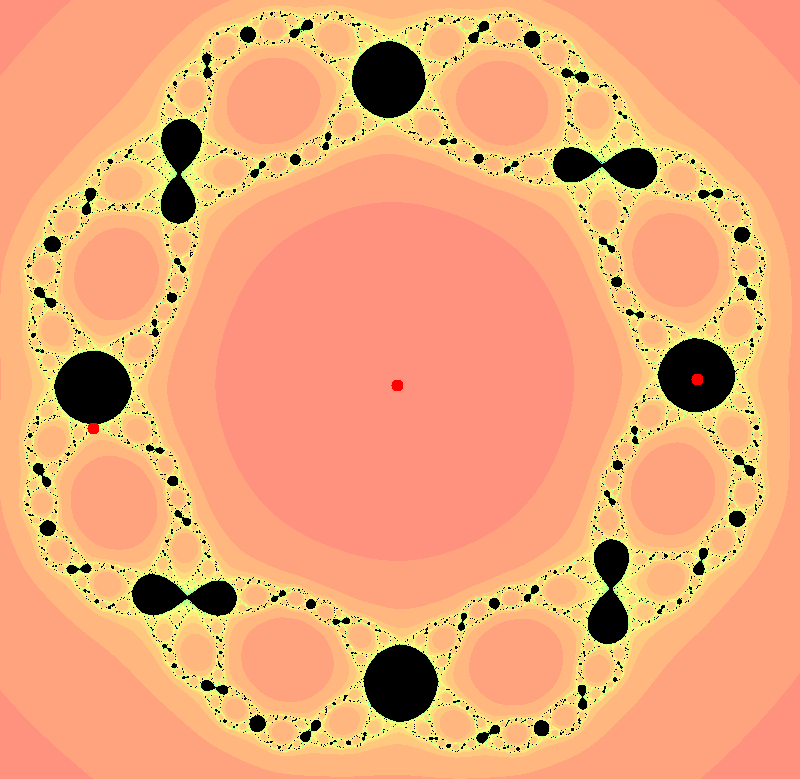}
\caption[Julia sets near bifurcation as $v_-$ leaves $K_-$]{Julia set for $\Rna,\,\, n=4, a=0.16+0.026i$ showing bifurcation as $v_-$ leaves $K_-$. $v_{\pm}, 0$ are marked in red, with $v_-$ on the left.}
\label{fig:lemniscatebifurcation}
\end{figure}

See Figure~\ref{fig:lemniscatebifurcation} for an example of this type, where the near figure 8 topological disks are alterations of the nearly perfectly round baby Julia sets (for quadratic maps from the main cardioid of the mandelbrot set).

\subsection{Paper Fortune Tellers: the basilica case}
\label{sec:PFT}

In this subsection, we recall the most general theorem from \cite{BoydBrouwer1} regarding angle re-identifications. 

To do so, we first need a bit of notation to help us refer to specific Fatou components contained in our baby basilica Julia sets and altered baby basilica Julia sets.

\begin{notation}
Given Assumptions~\ref{ass:general}, moreover assume $P_c$ is from the basilica bulb so that $c$ tends to an attracting 2-cycle.
For any connected portion of the Julia set which has been assigned external angles, we label a Fatou component $\FC$ by 
$\FClabel{a_1}{b_1}{a_2}{b_2}$
where $a_1 \sim b_2$ and $b_1 \sim a_2$ are the external angles that share landing points on the boundary of $\FC$ with lowest denominator, and where $0 \leq a_1 < b_1 < a_2 < b_2 \leq 1$. 
\end{notation}

\begin{figure}
\includegraphics[width=0.8\textwidth]{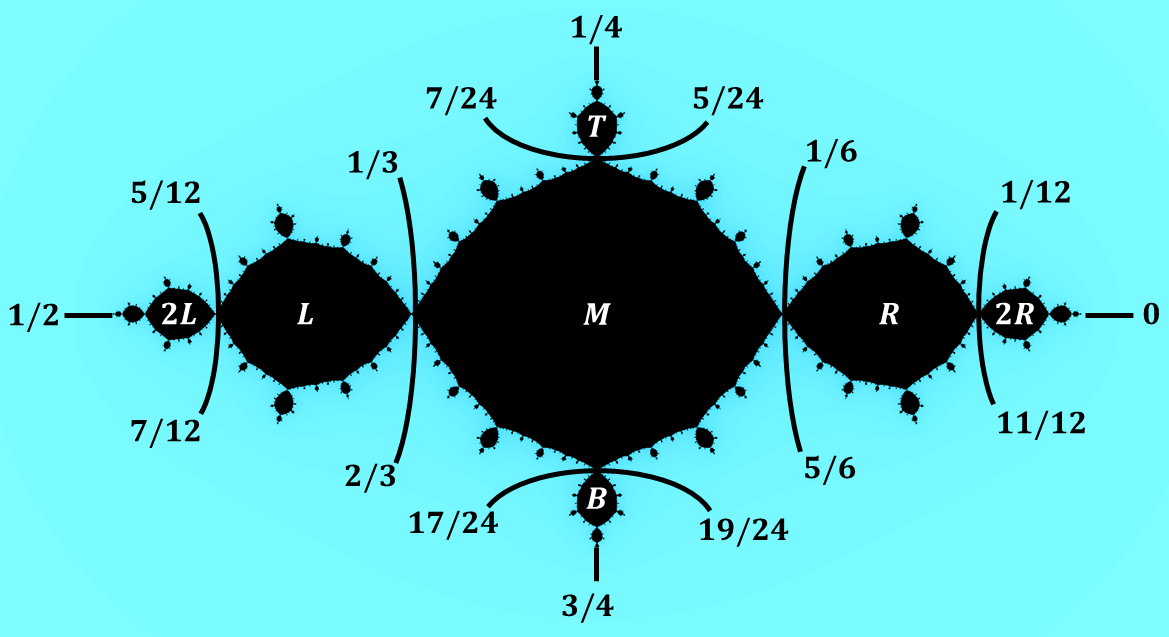}
    \includegraphics[width=0.8\textwidth]{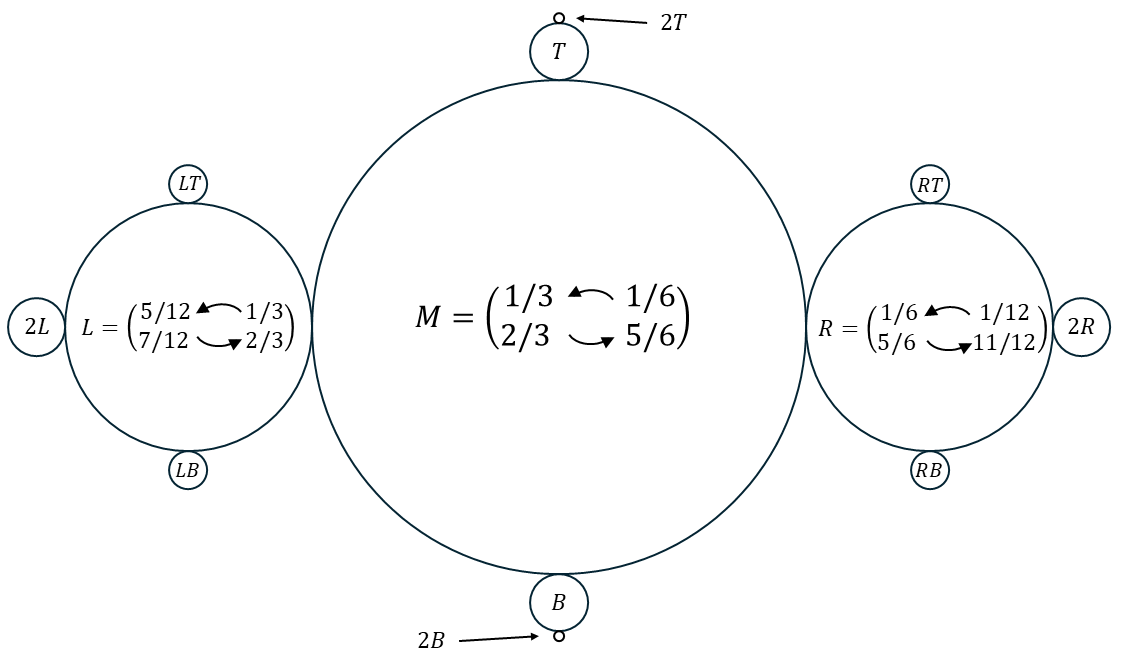}
    \caption{Standard basilica rays and diagram naming some Fatou components of the basilica Julia set in terms of rays.}
    \label{fig:labels}
\end{figure}

See Figure~\ref{fig:labels} for some selected ray assignments for the basilica and names of some major components. 
For Fatou components lying along the real axis, the point with $b_1 \sim a_2$ is directly left of the point with $a_1\sim b_2$ (the latter of which is the first point of the component encountered when walking from the point assigned angle $0$ counterclockwise around the Julia set). Then we think of the component as horizontally aligned and the two identifications can be considered ``vertical pinches'' separating this Fatou component from others nearby. 

Or, if a Fatou component bulb is not along the main horizontal axis but in the top half, i.e., is identified by angles all of which fall in $(0,1/2)$, then tilting one's head right will line up the bulb with the matrix label.
For any bulb not along the main line but with angles in $(1/2, 1)$, one needs to tilt one's head left to get the bulb to match this matrix notation. 

We refer to two Fatou components of the filled Julia set as \textbf{adjacent} if their boundaries touch.     
If two Fatou components share a boundary point, that point is the location of an angle identification. A simple example is that $M=\FClabel{1/6}{1/3}{2/3}{5/6}$ is adjacent to $L=
\FClabel{1/3}{5/12}{7/12}{2/3}
$ because the two components share the $\alpha$-fixed point as a boundary point, which is the point identified with both the angles ${1}/{3}$ and ${2}/{3}$. Typically, however, only one of two adjacent components will contain the identified angles of their shared boundary point in its constructed name; generally, the ``smaller'' component will. 

Now, in the case of the baby quadratic Julia set being a basilica, referring to Figure~\ref{fig:labels}, we use $M$ to refer to the central component: $M=\FClabel{1/6}{1/3}{2/3}{5/6}$,
and note $M$ is the location of the critical point in the basilica. We use $L$ and $R$ for the largest components adjacent to $M$, so 
$R=
\FClabel{1/12}{1/6}{5/6}{11/12}
$ and $L=
\FClabel{1/3}{5/12}{7/12}{2/3},
$ 
and $L$ is the location of the critical value in the basilica.  

Now, we can finally state a version of the most general theorem from \cite{BoydBrouwer1}. The following is slightly modified, since in \cite{BoydBrouwer1}, the theorem is stated in terms of lamination diagrams.

\begin{theorem} \label{thm:typeN} (\cite{BoydBrouwer1})
    Under Assumptions~\ref{ass:general}, assume $P_c$ is from the basilica bulb and  
    suppose that $v_-$ lies in a Fatou component $\FC_-$ of $\Kminus$ which corresponds to the Fatou component $\FC_c$ of $K_c$ such that $\FC_c$ is \textbf{not}
    $L=\FClabel{1/3}{5/12}{7/12}{2/3}$. If instead we have $v_-\in L$, then the preimage is not altered. 

    Starting with $\FC_0=L$ and ending with $\FC_N=\FC_c$, let $\FC_0, \dots, \FC_N$ be the shortest possible path of adjacent Fatou components in $K_c$ from $L$ to $\FC_c$. For $i=1, \dots, N$, let $\FC_i$'s preimage Fatou components be labeled $\FClabel{a_1^i}{b_1^i}{a_2^i}{b_2^i}$ and $\FClabel{a_3^i}{b_3^i}{a_4^i}{b_4^i}$ , noting that $a_1^i \sim b_2^i, b_1^i \sim a_2^i, a_3^i \sim b_4^i, b_3^i \sim a_4^i$. 

    Then the angle identifications for $\gamma_0$ are the same as the angle identifications for $\gamma_-$
    except that for each $i=1, \dots, N$, the identifications change to $a_1^i \sim b_4^i$ and $b_2^i \sim a_3^i$ (with the identifications $b_1^i \sim a_2^i, b_3^i \sim a_4^i$ unchanged), with one exception: if $\FC_1 = M=
    \FClabel{1/6}{1/3}{2/3}{5/6}$, 
    then $\gamma_0$ has the identified angles $b_1^1=\frac{1}{6} \sim \frac{1}{3}=a_3^1$ and $b_4^1=\frac{2}{3} \sim \frac{5}{6}=a_2^1$.
    \end{theorem}

This sequence of ray identification changes is effectively a sequence of paper fortune teller moves, along the path which the critical value takes to get from its ``expected'' location ($L$) to its actual location ($\FC_-$).

\section{Sidecars: $\nu$-Rabbits}
\label{sec:sidecars}

In this section, we examine identified angle changes that occur when a baby Julia set in the form of a ``rabbit'' is altered. The rabbit bulbs are the {primary decorations} of the main cardioid of the Mandelbrot Set, $\mandel$, by which we mean that the rabbit bulbs are hyperbolic components of $\mandel$ that are adjacent to the main cardioid. More specifically, we call a Julia set or quadratic map a \textbf{$\nu$-rabbit} when it associated with a parameter from a bulb attached to the main cardioid and has an attracting cycle of length $\nu$. 
The right side of Figure~\ref{fig: 3-rabbit v_- in TL pic} is a baby 3-rabbit.

Specifically, the case we study in this section satisfies the following. 

\begin{assumptions}
    \label{ass:rabbit type N} Let $\mapR=\Rnab$ be a generalized McMullen map satisfying Assumptions~\ref{ass:general}; and moreover, in (A4), $P_c$  is chosen from a $\nu$-rabbit bulb attached to the main cardioid of the Mandelbrot set $\mandel$, where $\nu\geq 3$ is an integer.

    We may abuse notation and refer to the corresponding components within $K(\mapRnab)$ by the same names.
\end{assumptions}

First, some helpful terminology. 

\begin{figure} \centering
     \includegraphics[height=6cm]{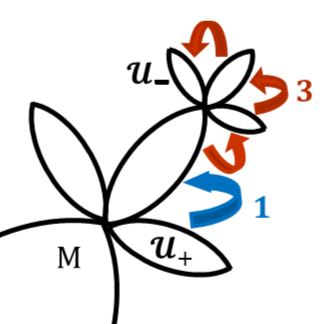}
    \caption[Counting rotations on path from $\FC_+$ to $\FC_-$]{Counting rotations on a sample path from a $\FC_+$ to a $\FC_2=\FC_-$ in a 4-rabbit. Here, $\FC_+=\FC_0$ is 1 rotation from $\FC_1$, and $\FC_1$ is 3 rotations from $\FC_-=\FC_2$. \label{fig:rabbit-U-path}}
\end{figure}

\begin{definition} \label{defn:nu-rabbit-terms}
 Under Assumptions~\ref{ass:rabbit type N}:
\begin{enumerate}
    \item[(a)] 
    For any $P_c$ taking the form of a $\nu$-rabbit, we refer to
    a point at which $\nu$ components of the filled Julia set meet as a \textbf{component junction}. 

   \item[(b)]
   A rabbit's \textbf{primary ears} are the largest Fatou components which meet the central component $M$ at a component junction. We may use \textbf{ear} for other non-main Fatou components as they are an ear of an ear in a sense. 
    
    \item[(c)] 
\textbf{Adjacent components} refer to any two Fatou components that meet at a component junction.

    \item[(d)] 
    As a Fatou component may have many adjacent components at each component junction, we may also note how many counterclockwise \textbf{rotations} one component is from another. (See Figure~\ref{fig:rabbit-U-path}.)

    \item[(e)] We take a \textbf{step} when we move from one Fatou component to an adjancent Fatou component.

\item[(f)] As in Notation~\ref{not:FC_- etc}, $v+$ is in a baby J-set $K_+$, and  $v_-$ is in a Fatou component $\FC_-$ of an eventual preimage $K_-$ of $K_+$.  

Additionally,  $\FC_+$ denotes the Fatou component of $K_-$ which is the ``expected'' location of the critical value; e.g., 
for the standard (upper) 3-rabbit, it is the small ear one counterclockwise rotation from $M$ on the upper side of the rabbit.
\end{enumerate}
\end{definition}

The number of $\nu$-rabbit bulbs for each $\nu\geq3$ is  $\varphi(\nu)$, where $\varphi$ is Euler's totient function. 
Using the map which takes the unit interval to the boundary of the Mandelbrot set $\mandel$, 
we  identify each rabbit bulb and, in fact, all hyperbolic components, by referencing the least of the external angles which land on the boundary of that component with the smallest denominators; for example, the standard upper 3-rabbit bulb is the $1/7$-bulb as the rays landing on the boundary of this bulb with smallest denominator are $1/7$ and $2/7$. 
Due to the way this mapping of $[0,1] \to \partial \cM$ was defined, we are able to use this angle again to identify the location of the critical value component $\FC_+$, which is always one of the rabbit's \text{ears}.

The reason  there are $\varphi(\nu)$ hyperbolic components which produce $\nu$-rabbits for each $\nu$ follows exactly from this construction. The attracting cycle for a $\nu$-rabbit is made up of $M$, $\FC_+$, and the other ears on the same ``side'' of the central component as $\FC_+$. The number of counterclockwise rotations between $M$, which contains the critical point, and $\FC_+$, which contains the critical value, sets a pattern that continues for the full length of the cycle. If $\FC_+$ were to be $k$ rotations away from $M$ in a $\nu$-rabbit where $k\vert \nu$, then not all $\nu$ components of the cycle would be visited. Thus the number of different ear positions that $\FC_+$ can hold is $\varphi(\nu)$, and each position corresponds to a different hyperbolic component of $\mandel$.

With this understanding of structure of quadratic rabbit Julia sets, we first provide our result that explains the alterations made to any baby $\nu$-rabbit preimage within $K(\Rnab)$. 

\begin{theorem} \label{thm: rabbit type N}

Let $\mapR=\Rnab$ be a generalized McMullen map satisfying the above Assumptions~\ref{ass:rabbit type N}.

    Starting with $\FC_0=\FC_+$ and ending with $\FC_N=\FC_-$, let $\FC_0, \dots, \FC_N$ be the shortest possible path of adjacent Fatou components in $\Kplus$ (or imagine it as $K(P_c)$) from $\FC_+$ to $\FC_-$. Then $N$ is the number of steps in the alteration. For each $i=1,\dots,N$ suppose  $\FC_{i}$ is $0\leq \ell_i \leq \nu-1$ counterclockwise rotations from $\FC_{i-1}$. See Figure~\ref{fig:rabbit-U-path}. 
    
    As we use letters to denote specific angles, we let 
    $\bA(x)$ denote the $x^{th}$ letter of the alphabet (if $\nu > 26$, you have to get creative and add letters to the alphabet).
       For $i=1,\dots,N$, each $\FC_i$ has two preimage components in $K(P_c)$ which share a common neighbor at the $i^{th}$ step. 
    Let $z_i^1$ and $z_i^2$ be the points at which these preimages meet the common neighbor, and let
       $a_i^1 \sim b_i^1 \sim \dots \sim \bA(\ell_{i})_i^1 \sim \bA(\ell_{i}+1)_i^1 \sim \dots \sim \bA(\nu)_i^1$ and $a_i^2 \sim b_i^2 \sim \dots \sim \bA(\ell_{i})_i^2 \sim \bA(\ell_{i}+1)_i^2 \sim \dots \sim \bA(\nu)_i^2$
    be the angles which are identified at each of these points, respectively, as they occur on $\Jminus$. 

    Then the identified angles on $\Jnot$ are the same as those on $\Jminus$, except that for each $i=1,\dots,N$, the first $\ell_i$ angles of each identification swap. Explicitly, we have that for each $i=1,\dots,N$, the
    identifications on $\Jnot$ change to
    $a_i^2 \sim b_i^2 \sim \dots \sim \bA(\ell_{i})_i^2 \sim \bA(\ell_{1}+1)_i^1 \sim \dots \sim \bA(\nu)_i^1$ and $a_i^1 \sim b_i^1 \sim \dots \sim \bA(\ell_{i})_i^1 \sim \bA(\ell_{i}+1)_i^2 \sim \dots \sim \bA(\nu)_i^2$
    if $\FC_{i}$ lies $\ell_i$ counterclockwise rotations from $\FC_{i-1}$, with all other identifications left unchanged.
\end{theorem}

To ease the reading of this notation, consider the subscript on component $\FC$, point $z$, and angle $\bA(x)$ to keep track of which step of the alteration each object corresponds to. Any object with a superscript is a preimage which occurs in $\Knot$ instead of $\Kminus$, where the number on the superscript is either a 1 or a 2 to correspond to the first or second preimage as they occur in $\Kminus$, as they will be altered to new groupings on $\Knot$. 

This result is best understood through examples, which we subsequently give organized primarily by the length of the attracting cycle and secondarily by the bulb of $\mandel$ related to the baby Julia set. For a simple example,
in an unaltered $\frac17$-bulb 3-rabbit as shown in Figure~\ref{fig: 3-rabbit sketch std upper}, the angles $1/7, 2/7, 4/7$ all meet at a point of the central component $M$, and the angles $9/14, 11/14,$ and $1/14 (= 15/14)$ all meet at a symmetrically opposite point of $M$. But, when $v_-$ is in the upper left ear of $K_-$, rather than upper right for the standard 3-rabbit, 
the angles $1/7, 1/14, 11/14$ meet at a point in the altered baby $J$, and the angles $2/7, 4/7, 9/14$ meet at a common point: see Figure~\ref{fig: 3-rabbit v_- in TL diagram} for a visual and Example~\ref{ex: 3-rabbit v_- in TL} for full details.

Without giving full proof, we note that this result is a straightforward adaptation of Theorem~\ref{thm:typeN}. The primary difference between this scenario and the one explained in Theorem~\ref{thm:typeN} is the presence of more than two angles at each identification. That being said, if one imagines the paper fortune teller motion that was used on the alteration of basilicas, this is the same concept, where two consecutive angles at each identification must split apart and reidentify, where there may be extra angles that travel with each side. For example, if five angles are identified at a point, the second and third angles would mark the start and end of the second component at that component junction, so splitting them would open that component up so it could be combined with another. This would be splitting the identified angles into groups of two and three which will stay together. Then to ensure that each identification maintains five angles in total, the groups of two must swap. Thus the extra angles simply stay with their part of the split, and the location of the split depends on which components needs to be opened up and combined to maintain the $2:1$ mapping of critical point onto critical value. The fact that we could write this down as always swapping the first $\ell$ angles had to do with the careful ordering of the angles.

\textbf{Consistency with basilicas: } 
In fact, Theorem~\ref{thm: rabbit type N} applies for all primary decorations of the main cardioid of $\mandel$, which includes the basilica 2-cycle bulb. The results for this bulb were phrased differently in \cite{BoydBrouwer1}, but when reworked to match the generalized $\nu$-cycle decorations, fall into the same patterns. In fact, in \cite{BoydBrouwer1}, the case in which $v_-$ passes through $M$ was set aside as distinct from other cases in terms of which angles reidentify; however, when arranging identified angles to always follow the rule $\frac{j}{3*2^k} \sim \frac{j+2}{3*2^k}$ where the second angle may be reduced mod $\mathbb{Z}$, this discrepancy is removed. (See Proposition~\ref{prop: 3-rabbit id angle sets}).


\subsection{3-rabbits}
We begin with the simplest rabbits: $3$-rabbits.  These  can be found in the $\frac17$-bulb and in the $\frac67$-bulb of the Mandelbrot set. These two cases are essentially symmetric reflections of each other, so we focus on $1/7$-rabbits in our examples below (for the interested reader, \cite{Brouwer-Thesis} contains several examples from the $6/7$-bulb, we highlight Example 5.17). But first, in the following subsection, we describe $3$-rabbits a bit more, before describing how they are altered in the subsequent subsection.

\subsubsection{Standard 3-rabbits}

First, we provide a well-known result about which angles can be identified together and how they relate to each other numerically. 

\begin{proposition} \label{prop: 3-rabbit id angle sets}
     Let $c$ be taken from the $\frac17$-bulb of the Mandelbrot set $\mathcal{M}$. Then $K(P_c)$ is a 3-rabbit, so at component junction $z$ within $K(P_c)$, the three angles that are identified at $z$ are either specifically $\frac17 \sim \frac27 \sim \frac47$, or the angles follow the rule of $\frac{j}{7*2^k} \sim \frac{j+2}{7*2^k} \sim \frac{j+6}{7*2^k}$ where $k\geq1$. Each fraction is irreducible, but in some cases, $\frac{j+6}{7*2^k}>1$, and hence will be written in its equivalent form mod $\ZZ$.

     Alternatively, let $c$ be taken from the $\frac67$-bulb of $\mandel$, which also produces 3-rabbits. At any component junction $z$ in $K(P_c)$, the three angles that are identified at $z$ are specifically $\frac37 \sim \frac57 \sim \frac67$ , or the angles follow the rule of $\frac{j}{7*2^k} \sim \frac{j+4}{7*2^k} \sim \frac{j+6}{7*2^k}$ where $k\geq1$, where each fraction is irreducible, and where fractions are listed in increasing order. Either all angles will be less than one, or $\frac{j+4}{7*2^k}$ and $\frac{j+6}{7*2^k}$ will both be greater than 1, in which case both fractions will be reduced mod $\ZZ$. 
\end{proposition}

As a brief example, on $\frac{1}{7}$-bulb Julia sets, we see that $\frac{9}{28} \sim \frac{11}{28} \sim \frac{15}{28}$, which matches the first pattern. These Julia sets also have $\frac{23}{28} \sim \frac{25}{28} \sim \frac{29}{28}$, but since this identification wraps around the point at which $0\equiv1$ lands, we write $\frac{1}{28}$ instead of $\frac{29}{28}$. Similarly, from the $\frac{6}{7}$-bulb, we would have $\frac{13}{14} \sim \frac{17}{14} \sim \frac{19}{14}$, but since the latter two fractions are greater than 1, we reduce them and write $\frac{13}{14} \sim \frac{3}{14} \sim \frac{5}{14}$.

In Section~\ref{sec:prelim}, we re-stated the convention of \cite{BoydBrouwer1} to give names to specific Fatou components of a baby Julia set or one of its preimages within the Julia set of $\Rnab$. However, the focus of that study was baby basilicas, which contain attracting 2-cycles and hence have angles that are identified in sets of 2. In general, for primary decorations of $\mandel$, a baby Julia set which contains an attracting $\nu$-cycle will have angles that are identified in sets of size $\nu$. While for these higher degree cycles we can still name a Fatou component by referring to two specific pairs of identified angles that lie on its boundary, we need to add more detail about which two of the three or more identified angles from each grouping should be chosen to name the component. Below, we give our naming convention for baby Julia sets taken from the 3-rabbit bulb. For $\nu>3$, we follow similar logic but do not state the identified angle names explicitly.

\begin{notation} \label{not: 3-rabbit comp labels}
    Let $\FC$ be a Fatou component of $K(P_c)$
    where $c$ is taken from the $\frac17$- or $\frac67$-bulb of the Mandelbrot set $\mathcal{M}$.
    Consider all points $z$ on $\partial \FC$ at which more than angle is identified. Let $a_1 \sim b_1 \sim c_1$ be the identified angles on $\partial \FC$ which have the largest possible denominator for all such $z$ and let $a_2 \sim b_2 \sim c_2$ have the second largest possible denominator. For each identification, assume that $a_j<b_j<c_j$. $\FC$ should be given an identified angle name as follows: 
    \begin{enumerate}
        \item If $a_1 < a_2 < c_2 < b_1$, the component should be labeled $\FClabel{a_1}{a_2}{c_2}{b_1}$.

        \item If $b_1 < a_2 < c_2 < c_1$, the component should be labeled $\FClabel{b_1}{a_2}{c_2}{c_1}$.
        
        \item If $c$ is from the $\frac17$-bulb of $\mandel$,
        we may have that $c_i$ has been reduced modulo $\mathbb{Z}$ for $i=1,2$. In this case, we have $c_2 < c_1 < b_1 < a_2$, and the component should be labeled $\FClabel{c_2}{c_1}{b_1}{a_2}$.

        \item On the other hand, if $c$ is from the $\frac67$-bulb of $\mandel$, we may have that $b_i$ and $c_i$ are both reduced modulo $\mathbb{Z}$ for $i=1,2$. In this case, we have $c_2 < b_1 < a_1 < a_2$, and the component should be labeled $\FClabel{c_2}{b_1}{a_1}{a_2}$.

        \item The central component containing the critical point should be labeled $M=\FClabel{1/14}{1/7}{4/7}{9/14}$ if $c$ is from the $\frac17$-bulb of $\mandel$ and  $M=\FClabel{5/14}{3/7}{6/7}{13/14}$ if $c$ is from the $\frac67$-bulb of $\mandel$, neither of which follows any of the previous naming conventions.
    \end{enumerate}
\end{notation}

We illustrate this proposition in the following example. 

\begin{figure} 
\begin{center}
    \includegraphics[height=9cm]{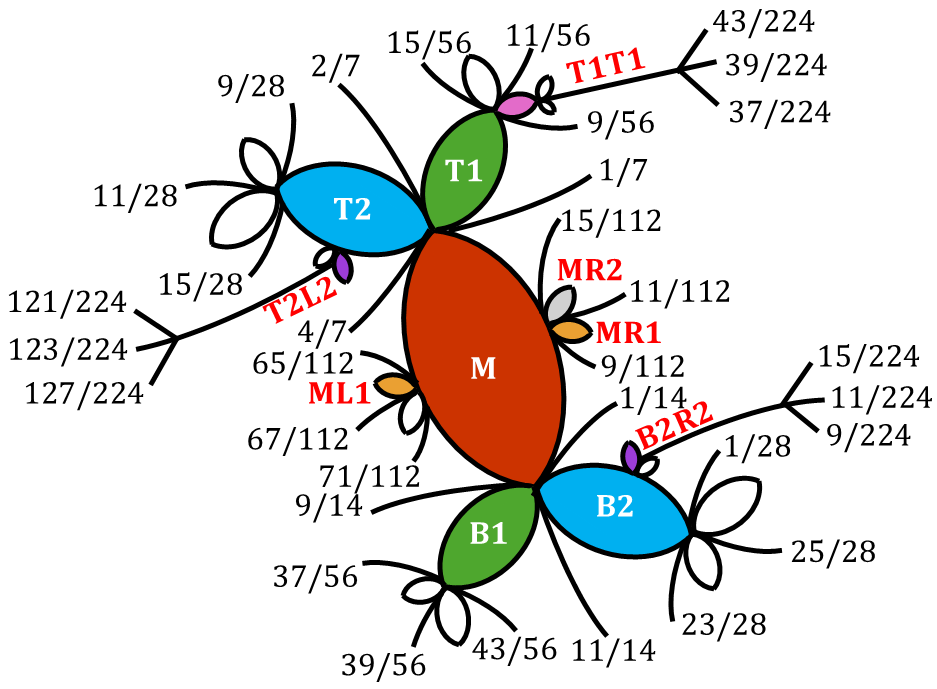}
\end{center}
\caption[External Angles standard (upper) 3-rabbit]{A diagram showing the identified angles on a 3-rabbit taken from the $\frac17$-bulb of $\mandel$\label{fig: 3-rabbit sketch std upper}.}
\end{figure}

\begin{example}
\label{ex: 3-rabbit comp labels}
    Consider any 3-rabbit Julia set taken from inside the $\frac17$-bulb of $\mandel$.

    As an example of the first case, the first ear encountered traveling counterclockwise from $M$ is $T1$, located in the ``top right'' position, which is the one that contains the critical value. On this component, the sets of identified angles which share the two lowest denominators are $\frac17 \sim \frac27 \sim \frac47$ and $\frac{9}{56} \sim \frac{11}{56} \sim \frac{15}{56}$. Here, we have that $\frac17 < \frac{9}{56} < \frac{15}{56} < \frac27$, so the component should be labeled $T1=\FClabel{1/7}{9/56}{15/56}{2/7}$.

    As an example of the second case, the second ear encountered when traveling counterclockwise from $M$ is the ``top left'' ear, $T2$, which has points with identified angles of largest denominator $\frac17 \sim \frac27 \sim \frac47$ and $\frac{9}{28} \sim \frac{11}{28} \sim \frac{15}{28}$. Since $\frac27 < \frac{9}{28} < \frac{15}{28} < \frac47$, the component should be labeled $T2=\FClabel{2/7}{9/28}{15/28}{4/7}$.

    As an example of the third case, the ``bottom right'' ear, $B2$ has points with identified angles of largest denominator $\frac{9}{14} \sim \frac{11}{14} \sim \frac{1}{14}$ and $\frac{23}{28} \sim \frac{25}{28} \sim \frac{1}{28}$. Notice that both of these sets of identified angles are arranged in increasing order following the description given in Proposition~\ref{prop: 3-rabbit id angle sets}, where each of the third angles have been reduced modulo $\mathbb{Z}$ (for example, $\frac{9+6}{14} = \frac{15}{14} \equiv \frac{1}{14} \text{ mod }\ZZ$). Arranging these angles in increasing order, we see that $\frac{1}{28} < \frac{1}{14} < \frac{11}{14} < \frac{23}{28}$, so the component should be labeled $B2=\FClabel{1/28}{1/14}{11/14}{23/28}$.

    To demonstrate the fourth case, we must instead take $c$ to be from a $\frac67$-bulb 3-rabbit. In this Julia set, the attracting cycle occurs between the central component and the two lower components, and the angle 0 lands on a component that stems off of the top right ear $T1$. The identified angles of largest two denominators on $T1$ are $\frac{13}{14} \sim \frac{3}{14} \sim \frac{5}{14}$ and $\frac{27}{28} \sim \frac{3}{28} \sim \frac{5}{28}$. Since $\frac{5}{28} < \frac{3}{14} < \frac{13}{14} < \frac{27}{28}$, the component should be labeled $\FClabel{5/28}{3/14}{13/14}{27/28}$.

    A diagram showing some identified angles on an unaltered 3-rabbit from the $\frac17$-bulb is provided in Figure~\ref{fig: 3-rabbit sketch std upper} for visual confirmation. 
\qed \end{example}

\subsubsection{Altered 3-rabbits} 
\label{subsub: 1/7 3-rabbits}

The examples in this sub-sub-section satisfy:

\begin{assumptions}\label{ass:17bulb}
    $\Rnab$ satisfies Assumptions~\ref{ass:general}, and further, $\Rnab$ is conjugate on $\Kplus$ to $P_c$ on $K(P_c)$ for a $c$ from the $\frac17$-bulb of 3-rabbits. 
\end{assumptions}

First, we give some general information to help navigate baby 3-rabbits that follow the pattern of those from the $\frac17$-bulb. In this bulb, the central component containing the critical point is $M=\FClabel{1/14}{1/7}{4/7}{9/14}$, and the component containing the critical value is $\FC_+=T1=\FClabel{1/7}{9/56}{15/56}{2/7}$. We may refer to $\FC_+$ as the ``expected'' location of the critical value, as our alterations occur when the other critical value $v_-$ is positioned elsewhere within its preimage $\Kminus$ of $\Kplus$. Visually, it is the ``upper'' ears of $\Kplus$ that contain the attracting cycle, and each component in the cycle is positioned one counterclockwise rotation from the previous one. 
In Figure~\ref{fig: 3-rabbit sketch std upper} given previously, $M$ is colored red and $\FC_+$ is colored green. The other ear in this cycle is $T2=\FClabel{2/7}{9/28}{15/28}{4/7}$, colored blue, so that in terms of colors, the three-cycle is red $\to$ green $\to$ blue. 
We chose to color components if they mapped onto the attracting cycle, following the order of which colors map where. For example, if a component $\FC$ of $\Jnot$ was found to map on the blue component of $\Jminus$, then $\FC$ was colored green. This helps identify when key components had their external angles altered or when they remained the same, and it makes it easier to keep track of key components when their position was altered significantly.

When we reference Theorem~\ref{thm: rabbit type N}, here we use that $\nu=3$, so components on the path from $\FC_+$ to $\FC_-$ can only be one or two rotations away from the previous one. If a component $\FC_i$ on this path is one rotation from the previous component $\FC_{i-1}$, the angles related to this step will exchange only the first angle in each set. Otherwise, a component $\FC_i$ must be two rotations from $\FC_{i-1}$, which means the first two angles in each identification will swap. Sometimes this latter change is more noticeable if one considers it as only the last angle from each identification swapping. 

\smallskip

We first give the simplest type of example: one step from $\FC_+$ to $\FC_-$, and that step is a single counterclockwise rotation.

\begin{figure}
    \centering
    \includegraphics[width=\textwidth]{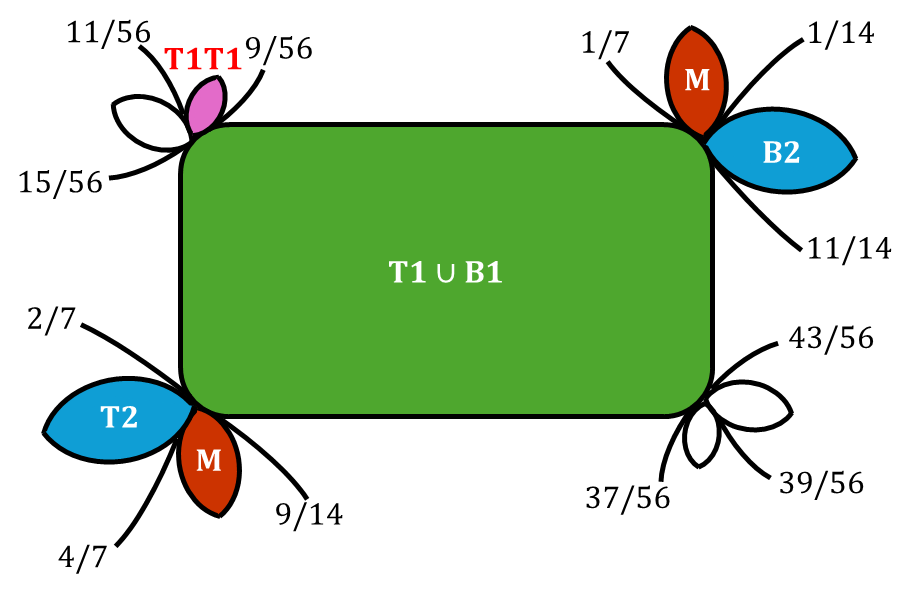}
    \caption[External Angles altered $\frac17$-bulb 3-rabbit, for $v_-\in TL$ ]{A diagram showing the external angle identifications on an altered $\frac17$-bulb 3-rabbit when $\FC_-=\FClabel{2/7}{9/28}{15/28}{4/7}$ \label{fig: 3-rabbit v_- in TL diagram}.}
\end{figure}

\begin{figure}
    \centering
    \begin{minipage}{0.45\textwidth}
        \includegraphics[width=\textwidth]{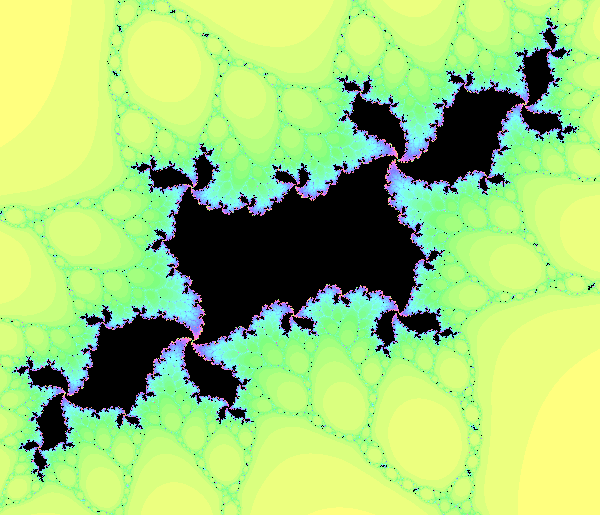}
    \end{minipage}
    \begin{minipage}{0.45\textwidth}
        \includegraphics[width=\textwidth]{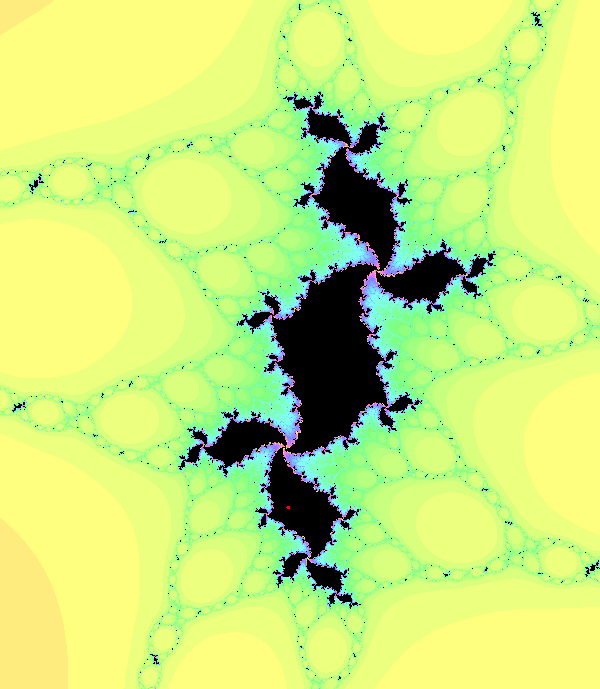}
    \end{minipage}
    \caption[Altered Julia for $\frac17$-bulb 3-rabbit, for $v_-\in T2$]{Both for $n=5$, $b=0.01+0.03i$, and $a=0.17353+0.04181i$, $\Kminus$ with $v_-$ identified inside of $T2$ (right), and its altered preimage $\Knot$ (left)\label{fig: 3-rabbit v_- in TL pic}.}
\end{figure}

\begin{example} \label{ex: 3-rabbit v_- in TL}
    Suppose $\Rnab$ satisfies Assumption~\ref{ass:17bulb}, and that $v_-$ is in a Fatou component $\FC_-$ of $\Kminus$ which we refer to as in the ``top left'' position, labeling it as $\FC_1=\FC_-=T2 = \FClabel{2/7}{9/28}{15/28}{4/7}$. This component is adjacent to the expected position of $v_-$, $\FC_+=T1$, as they share a common boundary point which is a component junction. Thus this example is of type $N=1$ since the components are one step apart. Furthermore, that step is a single counterclockwise rotation, which coincidentally mimics the behavior of the attracting periodic cycle.
    
    The Fatou component $\FC_1=T2$ of $K_c$ has two preimages within $K_c$, which are $\FC_1^1=\FClabel{1/7}{9/56}{15/56}{2/7}$ and $\FC_1^2=\FClabel{9/14}{37/56}{43/56}{11/14}$. These preimage components are both adjacent to $M
    $, and meet $M$ at the points $z_1^1$ and $z_1^2$, which are the locations of the identified angles $\frac17 
    \sim \frac27 \sim \frac47$ and $\frac{9}{14} \sim \frac{11}{14} \sim \frac{1}{14}$, respectively $a_1^1\sim b_1^1 \sim c_1^1$ and $a_1^2 \sim b_1^2 \sim c_1^2$ in the notation. One can see these components colored in green on Figure~\ref{fig: 3-rabbit sketch std upper}. Our goal is to combine these two preimages into one, so that the critical value component has a single preimage which maps onto it 2:1. 
    Following Theorem~\ref{thm: rabbit type N}, the fact that $\FC_-$ is a single rotation from $\FC_+$ means that only the first angle from each identification should swap, and only one set of angles will change identifications.
    This indeed meets that criteria, as on $\Jnot$
    we have $\frac{9}{14} \sim \frac27 \sim \frac47$ and $\frac17 \sim \frac{11}{14} \sim \frac{1}{14}$, with all other identifications from $\Jminus$ left unchanged. What we see is that this executes our goal of combining $\FC_1^1$ and $\FC_1^2$ into a single component, splitting their common neighbor of $M$ apart into two new components. This is the paper fortune technique carrying over from \cite{BoydBrouwer1}. What has changed is that there are now components which meet at the same component junctions as $\FC_1^1$ and $M$ or $\FC_1^2$ and $M$, but are not involved in this splitting and rejoining. We see that these components appear to stay attached at these component junctions, keeping their original identified angle names and relative positions in comparison to the components that split and combine. It is from this visual that we dub these components ``\textbf{sidecars}'', moving along with the changes and simply moving to a new position in order to keep the external angles arranged in increasing order around the full boundary of $\Jnot$. 
    
    A diagram showing the identified angles and a set of computer-generated images of such a $\Kminus$ and its corresponding $\Knot$ can be found in Figures~\ref{fig: 3-rabbit v_- in TL diagram} and~\ref{fig: 3-rabbit v_- in TL pic}, respectively. In the diagram, the combined new component is colored green to show that it is the green components from Figure~\ref{fig: 3-rabbit sketch std upper} which have combined, and the new components that split are colored red since they were formerly the single red component $M$. The blue components are the sidecars, maintaining their identified angle names from $\Jminus$ and simply moving as the component junctions move. 
\qed \end{example}

\begin{definition} \label{defn:sidecar}
    A \textbf{sidecar} is a Fatou component in an altered preimage of a baby Julia set, which is adjacent at the component junction to the \textbf{altered} components which split and recombine, but which keeps its original identified angle name and relative position.
\end{definition}

In the next example, we still have $\FC_-$ one step from $\FC_+$, moving between components at the same component junction as in Example~\ref{ex: 3-rabbit v_- in TL}, but now $\FC_-$ is two rotations away, which is the middle component $M$.

\begin{figure}
    \centering
    \includegraphics[height=10cm]{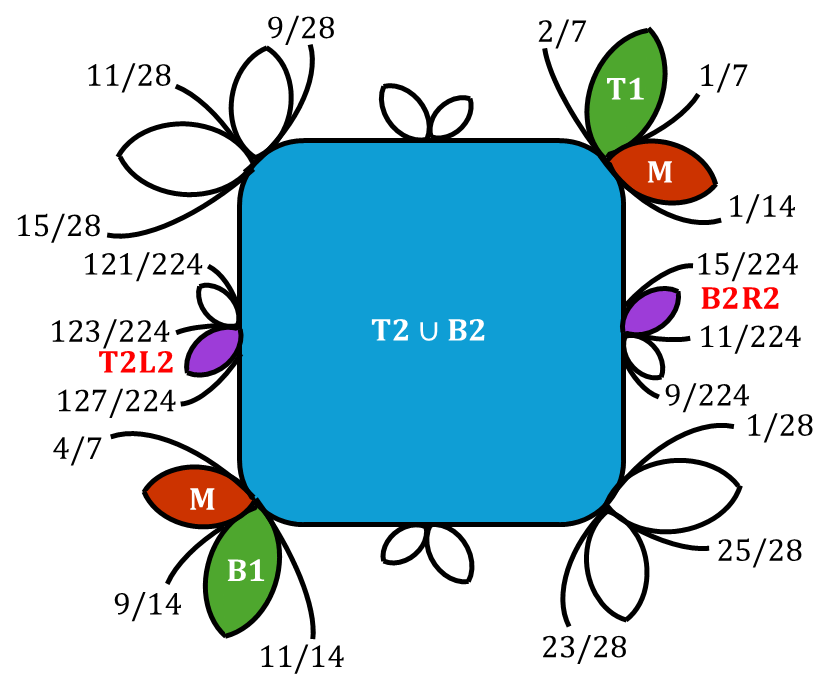}
    \caption[External Angles for altered 3-rabbit from the $\frac17$-bulb when $\FC_-=M$]{A diagram showing the external angle identifications present on an altered 3-rabbit from the $\frac17$-bulb when $\FC_-=M=\FClabel{1/14}{1/7}{4/7}{9/14}$\label{fig: 3-rabbit v_- in M diagram}.}
\end{figure}

\begin{figure}
    \centering
    \begin{minipage}{0.45\textwidth}
        \includegraphics[width=\textwidth]{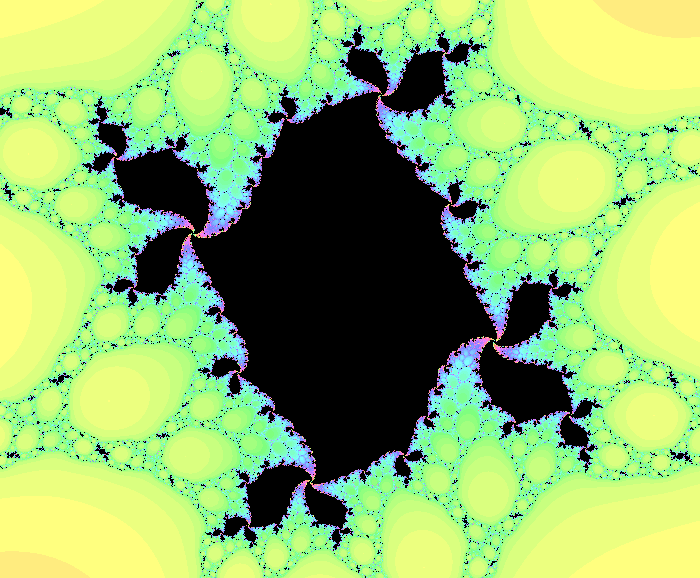}
    \end{minipage}
    \begin{minipage}{0.45\textwidth}
        \includegraphics[width=\textwidth]{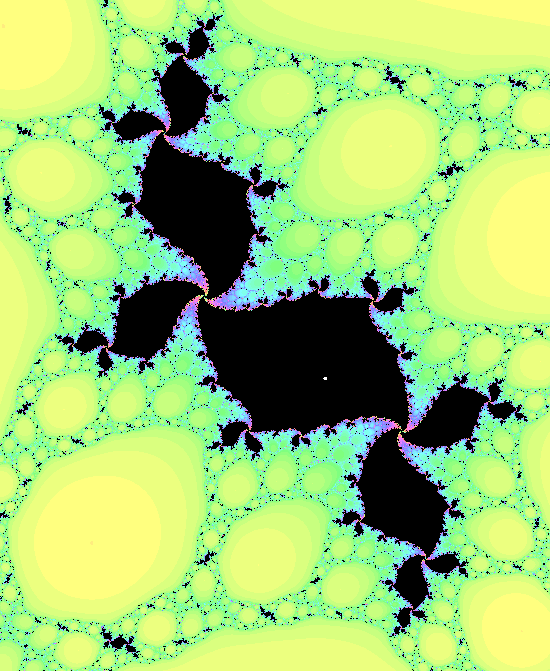}
    \end{minipage}
\caption[Altered Julia 3-rabbit from the $\frac17$-bulb when $\FC_-=M$]{\label{fig: 3-rabbit v_- in M pic}Both for $n=3$, $b=0.01+0.04i$, and $a=0.1-0.76i$, $\Kminus$ with $v_-$ identified inside of $M$ (right), and its altered preimage $\Knot$ (left).} 
\end{figure}

\begin{example} \label{ex: 3-rabbit upper v_- in M}
    Suppose $\Rnab$ satisfies Assumption~\ref{ass:17bulb}, and that
    $v_-$ is in the central component of $\Jminus$, which we refer to $M$ as in ``middle'' or ``main''; that is, let $\FC_-=M=\FClabel{1/14}{1/7}{4/7}{9/14}$. This component is also adjacent to $\FC_+=T1$, sharing the same common boundary point as $T1$ shares with $T2$, so this example is of Type $N=1$. However, in this case, $\FC_-=\FC_1$ is two rotations away from $\FC_+$.
    
    The two preimages of $M$ are $\FC_1^1=\FClabel{2/7}{9/28}{15/28}{4/7}$ and $\FC_1^2=\FClabel{11/14}{23/28}{1/28}{1/14}$, and as in Example~\ref{ex: 3-rabbit v_- in TL}, they are both adjacent to $M$, in fact using the same shared boundary points as in that example. Thus, the same sets of identified angles will be split and reidentified, but in different pairings since $M$ is two rotations from $T1$. In this example, the reidentified sets are $\frac{9}{14} \sim \frac{11}{14} \sim \frac47$ and $\frac17 \sim \frac27 \sim \frac{1}{14}$, which matches Theorem~\ref{thm: rabbit type N} in that the first two angles from each identification swap. 
    A diagram showing the identified angles and a computer-generated image of $\Knot$ beside the unaltered $\Kminus$ can be found in Figures~\ref{fig: 3-rabbit v_- in M diagram} and~\ref{fig: 3-rabbit v_- in M pic}, respectively. In Figure~\ref{fig: 3-rabbit v_- in M diagram}, we see that the new central component is colored blue, a combining of the two blue components from Figure~\ref{fig: 3-rabbit sketch std upper}, and the two new red components are the result of $M$ splitting apart.
\qed \end{example}

The next example documents that even when movement of one step away from $\FC_+$ occurs to an adjacent component that is not in the periodic cycle of Fatou components (i.e., not a ``primary ear''), we still get the same rule about counting the number of counterclockwise rotations.

\begin{figure}
    \begin{minipage}{0.45\textwidth}
        \includegraphics[width=\textwidth]{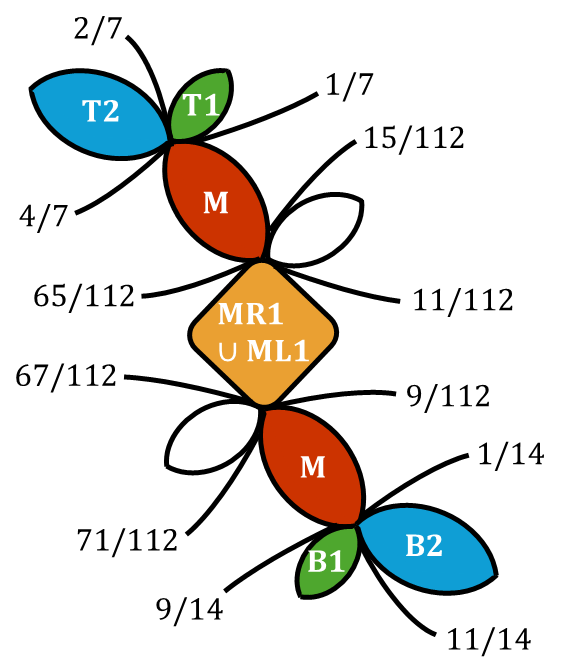}
    \end{minipage}
    \begin{minipage}{0.45\textwidth}
        \includegraphics[width=\textwidth]{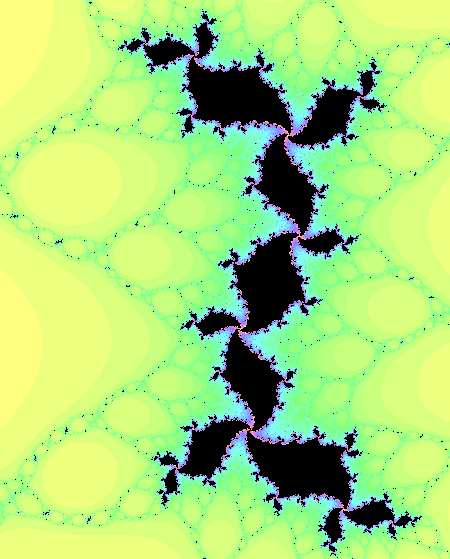}
    \end{minipage}
    \caption[External Angles and Julia for altered 3-rabbit from $\frac17$-bulb where $\FC_-=T1T1$.]{Left: A diagram showing the identified angles on an altered 3-rabbit taken from the $\frac17$-bulb where $\FC_-=T1T1=\FClabel{9/56}{37/224}{43/224}{11/56}$. Right: An example of left; specifically, $\Knot$ where $n=5$, $b=0.01+0.03i$, and $a=0.17357+0.041914i$.\label{fig: 3-rabbit upper v_- in (TR)(TR) altered}}
\end{figure}

\begin{example}\label{ex: 3-rabbit upper v_- in (TR)(TR)}
    Under Assumption~\ref{ass:17bulb}, suppose $\FC_-=\FClabel{9/56}{37/224}{43/224}{11/56}$,
    which we call $T1T1$, as this component is the first counter-clockwise ear* of the largest pair coming off of $T1$, where the pair is located at the top of $T1$ when $M$ is positioned in the center. $T1T1$ is colored pink in Figure~\ref{fig: 3-rabbit sketch std upper}.
    
    This component is one rotation away from $T1$, so this example is also of Type $N=1$ with a single rotation. The preimages of $\FC_-$ in $K(P_c)$ meet $M$ at the points with identified angles $\frac{9}{112} \sim \frac{11}{112} \sim \frac{15}{112}$ and $\frac{65}{112} \sim \frac{67}{112} \sim \frac{71}{112}$ and are colored orange in Figure~\ref{fig: 3-rabbit sketch std upper}, labeled as $MR1$ and $ML1$ (read as ``main right 1'' and ``main left 1''). As expected, we see that only the first angle in each of these identified pairs needs to change identification to combine the two preimage components into one. Explicitly, we see that on $\Jnot$, the identified angles previously listed have changed to $\frac{65}{112} \sim \frac{11}{112} \sim \frac{15}{112}$ and $\frac{9}{112} \sim \frac{67}{112} \sim \frac{71}{112}$, while all other identifications remain unchanged. A diagram of the external angles on $\Knot$ and a computer generated image of $\Knot$ are provided in Figure~\ref{fig: 3-rabbit upper v_- in (TR)(TR) altered}. Note that the central component on $\Knot$ has been colored orange, as the two orange components from Figure~\ref{fig: 3-rabbit sketch std upper} have combined, splitting $M$ apart into two red components. Since this paper fortune teller movement runs across the middle of $M$, we see that the blue and green ears are simply pushed out farther from the center unaltered, and that the components next to $MR1$ and $ML1$ from Figure~\ref{fig: 3-rabbit sketch std upper} are the sidecars in this case.
\qed \end{example}

We conclude this section with an example of Type $N=2$. Such examples were more difficult to complete, as the relevant identified angles got very small very quickly, but we include one to illustrate our result for a greater number of steps.

\begin{figure}
    \begin{minipage}{0.45\textwidth}
        \includegraphics[width=\textwidth]{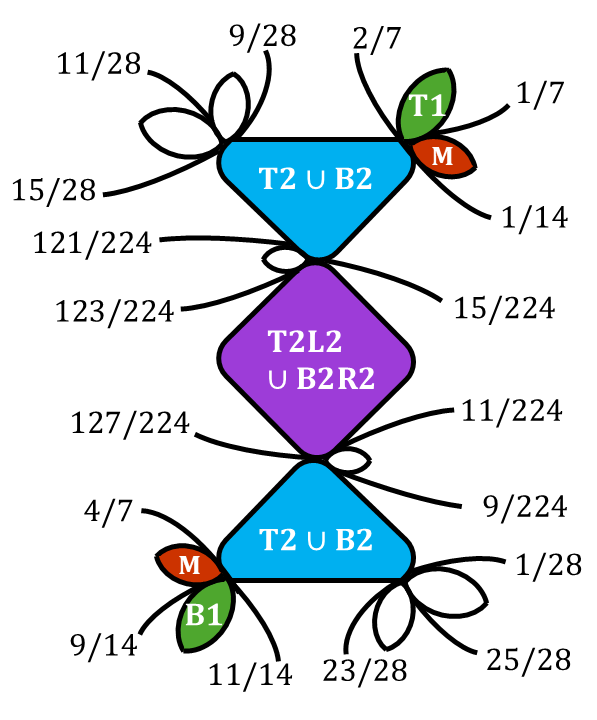}
    \end{minipage}
    \begin{minipage}{0.45\textwidth}
        \includegraphics[width=\textwidth]{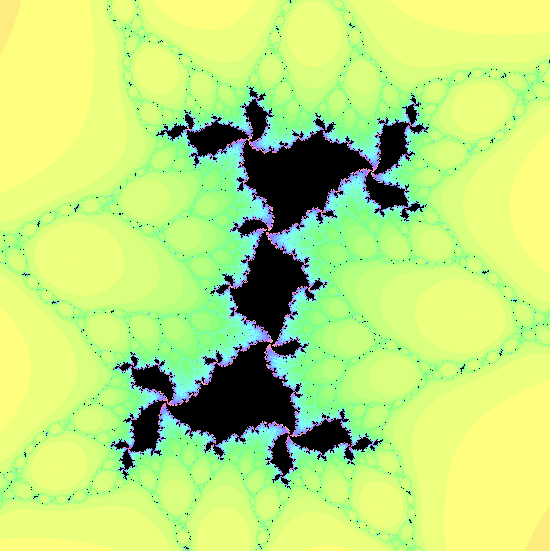}
    \end{minipage}
    \caption[External Angles and Altered Julia for upper bulb 3-rabbit with $\FC_-=MR1$]{Left: A diagram showing the identified angles on an altered 3-rabbit taken from the $\frac17$-bulb where $\FC_-=MR1=\FClabel{11/112}{23/224}{29/224}{15/112}$. Right: An example of left; specifically, $\Knot$ where $n=5$, $b=0.01+0.03i$, and $a=0.17356+0.041864i$.\label{fig: 3-rabbit upper v_- in M(RT) altered}}
\end{figure}

\begin{example}\label{ex: 3-rabbit upper v_- in M(RT)}
    Suppose that $\Rnab$ satisfies Assumption~\ref{ass:17bulb}, and that $\FC_-$ is the upper of the largest set of ears* on the right side of $M$, that is $\FC_-=MR2=\FClabel{11/112}{23/224}{29/224}{15/112}$, colored gray in Figure~\ref{fig: 3-rabbit sketch std upper}. 
    Note that this component is two steps away from $\FC_+=T1$ and is two rotations away at the second step. The preimage components of $MR2$ do not share a common neighbor on $\Kminus$, but they do share a common neighbor if the identification changes detailed in Example~\ref{ex: 3-rabbit upper v_- in M} are made, taking the movement of $v_-$ from $T1$ to $M$ to be the first step. Then the preimage components of $MR2$ meet their common neighbor at the identified angles $\frac{9}{224} \sim \frac{11}{224} \sim \frac{15}{224}$ and $\frac{121}{224} \sim \frac{123}{224} \sim \frac{127}{224}$. Since $\FC_2=MR2$ is two rotations away from $\FC_1=M$, we find that the first two angles in each identified set swap. Therefore, these angles are reidentified in $\Jnot$ as $\frac{121}{224} \sim \frac{123}{224} \sim \frac{15}{224}$ and $\frac{9}{224} \sim \frac{11}{224} \sim \frac{127}{224}$. A diagram detailing these angle identification changes and a computer generated image of $\Knot$ are given in Figure~\ref{fig: 3-rabbit upper v_- in M(RT) altered}. For this figure, track the purple components  (B2R2 and T2L2)
    from Figures~\ref{fig: 3-rabbit sketch std upper} and~\ref{fig: 3-rabbit v_- in M diagram}.
\qed \end{example}

\subsection{5-rabbits}

We note \cite{Brouwer-Thesis} contains several examples of $4$-rabbits. However, this case is functionally equivalent to the 3-rabbit case, just that there are two sidecars instead of one,  so we omit these examples here. 

In this subsection, we provide a couple of 5-rabbit examples demonstrating Theorem~\ref{thm: rabbit type N}. Since $\varphi(5)=4$, there are 4 different bulbs from which 5-rabbits are generated. The 5-rabbits in the 1/31-bulb act like the 3-rabbits in the 1/7-bulb, just with an extra couple of sidecars. So, we examine the 9/31-bulb, where the periodic orbit cycles through the 5 connected ears in an ``every other'' pattern; that is, if you follow the attracting periodic cycle, each Fatou component is two counterclockwise rotations from the previous one. We will not discuss 5-rabbits from the two lower bulbs, as each is functionally a reflection of its corresponding upper bulb.

\begin{assumptions}
    \label{ass:9-31-rabs}
$\Rnab$ satisfies Assumptions~\ref{ass:general}, and is conjugate on $K_+$ to $P_c$ with $c$ from the $9/31$-bulb of $\cM$.
\end{assumptions}

\begin{figure}
    \includegraphics[width=\textwidth]{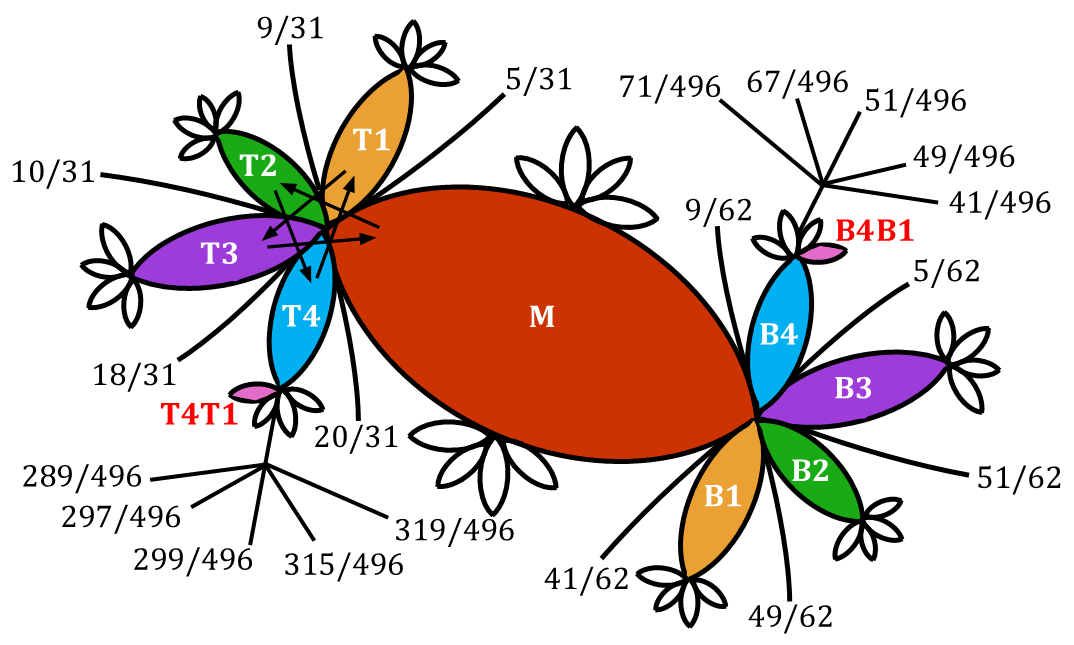}
    \caption[External Angles on the 5-rabbit from the $9/31$-bulb]{A diagram showing the angle identifications on a Julia set taken from the $\frac{9}{31}$-bulb of $\mandel$.\label{fig: std 9/31}}
\end{figure}

A diagram showing the external angle identifications on an unaltered 9/31-bulb 5-rabbit is given in Figure~\ref{fig: std 9/31}. All 5-rabbits taken from the $\frac{9}{31}$-bulb have $M=\FClabel{9/62}{5/31}{20/31}{41/62}$ and $\FC_+=\FClabel{9/31}{289/992}{319/992}{10/31}$, which we call $T2$ for its position two rotations away from $M$ located on the top set of largest ears.
The pattern of colors we use to represent the sequence of the components in the attracting cycle is red $\to$ green $\to$ blue $\to$ orange $\to$ purple. As mentioned above, for 9/31-bulb rabbits, these components are not arranged in order but in an alternating pattern, but we name the components in increasing order to keep their relative position clear over their placement in the cycle.

First, we consider an example of the identified angle changes that occur when $\FC_-$ is a component of the attracting cycle other than $\FC_+$, so that $\FC_-$ is one step away, but that step is four rotations. For such examples, we use the fact that the ears which comprise the attracting cycle and their symmetric components meet $M$ at the points where the angles $\frac{5}{31} \sim \frac{9}{31} \sim \frac{10}{31} \sim \frac{18}{31} \sim \frac{20}{31}$ and $\frac{41}{62} \sim \frac{49}{62} \sim \frac{51}{62} \sim \frac{5}{62} \sim \frac{9}{62}$ are identified.

\begin{example} \label{ex: 9/31 v_- in TR}
    Suppose $\Rnab$ satisfies Assumptions~\ref{ass:9-31-rabs}, and suppose $\FC_-=T1=\FClabel{5/31}{41/248}{71/248}{9/31}$, the upper of the orange ears in Figure~\ref{fig: std 9/31}. Thus, $v_-$ is four counterclockwise rotations away from the standard critical value location of $\FC_+=T2$, the upper of the green ears in Figure~\ref{fig: std 9/31}.
The preimage Fatou components of this orange
ear are the two blue ones, located between $\frac{18}{31}$ and $\frac{20}{31}$ and between $\frac{9}{62}$ and $\frac{5}{62}$, so it is these angle pair identifications which must split and re-combine to form one Fatou component. Swapping the first four angles from each identification, we have that
the identified angles on $\Jnot$ are $\frac{41}{62} \sim \frac{49}{62} \sim \frac{51}{62} \sim \frac{5}{62} \sim \frac{20}{31}$ and $\frac{5}{31} \sim \frac{9}{31} \sim \frac{10}{31} \sim \frac{18}{31} \sim \frac{9}{62}$.  
    See Figure~\ref{fig: 9/31 (TR)} for a diagram and computer-generated image of such a $\Jnot$.
    \qed

    \begin{figure} \centering
            \includegraphics[width=.8\textwidth]{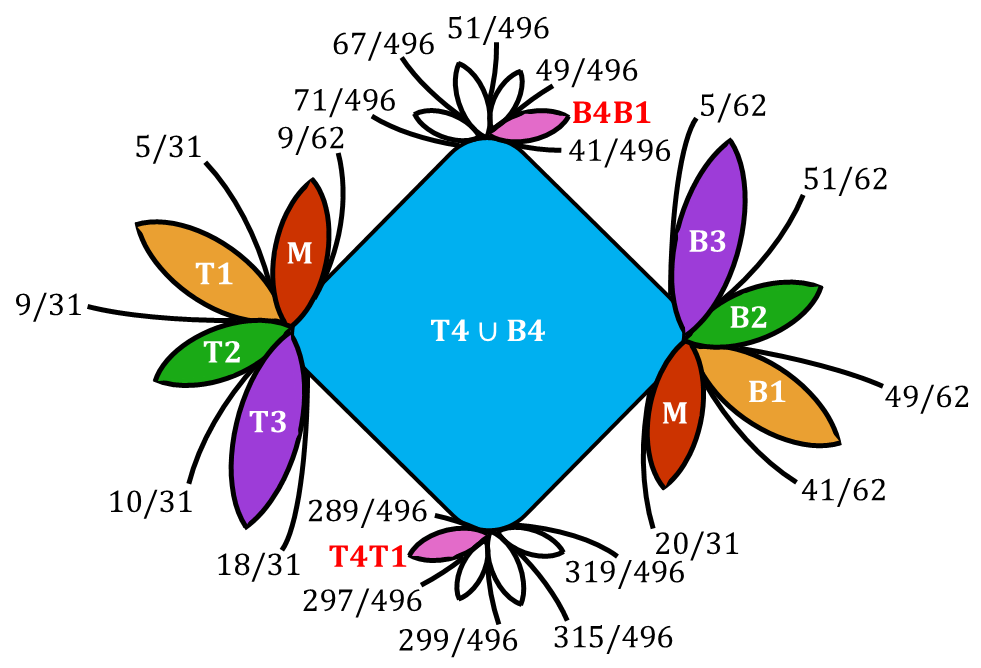}
            \includegraphics[width=.6\textwidth]{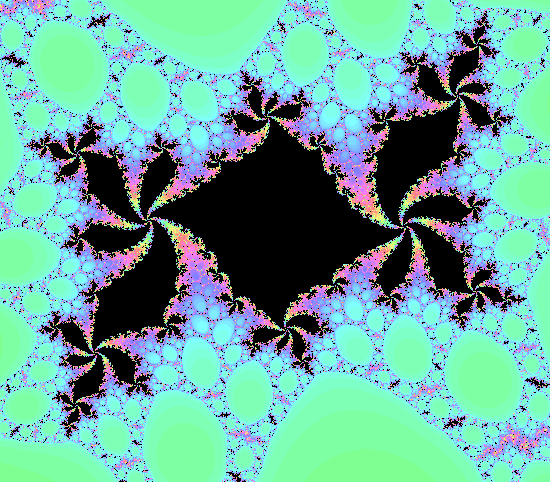}
        \caption[External Angles and Altered Julia on 5-rabbit from $\frac{9}{31}$-bulb where $\FC_-=T1$]{Top: A diagram of an altered 9/31-bulb 5-rabbit where $\FC_-=T1$. Bottom: $\Knot$ where $n=3$, $b=0.03+0.02i$, and $a=0.07008+0.02525i$.\label{fig: 9/31 (TR)}}
    \end{figure}
 \end{example}

Our second and final $5$-rabbit example builds on the prior example and is of Type $N=2$.

\begin{example} \label{ex: 9/31 N=2}
    Suppose $\Rnab$ satisfies Assumptions~\ref{ass:9-31-rabs}, and $\FC_-=T1T1=\FClabel{41/248}{165/992}{195/992}{49/248}$, which is the first counterclockwise rotation ear* of the largest set extending off of $T1$.
    This $\FC_-$ is two steps away from $\FC_+$, where $\FC_1=T1$ is four rotations from $\FC_+=\FC_0$ and $\FC_-=\FC_2$ is one rotation from $\FC_1$.  The preimages of $\FC_2=T1T1$ within $K(P_c)$ are colored pink in Figure~\ref{fig: std 9/31}, which share a common neighbor in Figure~\ref{fig: 9/31 (TR)} following the changes made in Example~\ref{ex: 9/31 v_- in TR} which we take as a first step. At this point, the angles at which these pink components (B4B1, T4T1) meet this common neighbor can now split and reidentify. Hence the angles identified on our final $\Jnot$ are the same as on $J(P_c)$ with the exception of the identifications $\frac{41}{62} \sim \frac{49}{62} \sim \frac{51}{62} \sim \frac{5}{62} \sim \frac{20}{31}$ and $\frac{5}{31} \sim \frac{9}{31} \sim \frac{10}{31} \sim \frac{18}{31} \sim \frac{9}{62}$ as a first step and $\frac{289}{496} \sim \frac{49}{496} \sim \frac{51}{496} \sim \frac{67}{496} \sim \frac{71}{496}$ and $\frac{41}{496} \sim \frac{297}{296} \sim \frac{299}{496} \sim \frac{315}{496} \sim \frac{319}{496}$ as a second step. These angle identifications changes are laid out in the diagram in Figure~\ref{fig: 9/31 (TR)(TR)}. \qed
    
    \begin{figure} \centering
            \includegraphics[width=0.9\textwidth]{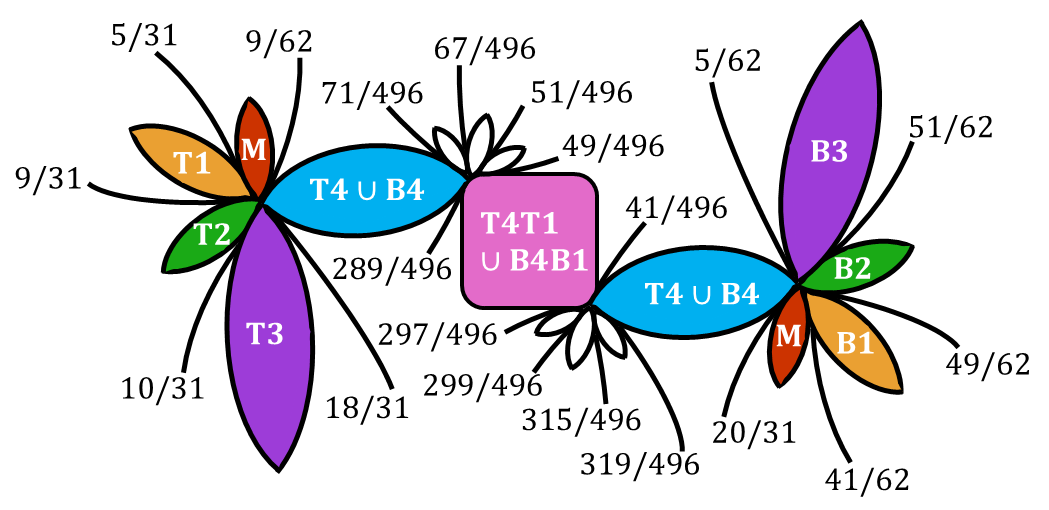}
            \includegraphics[width=.9\textwidth]{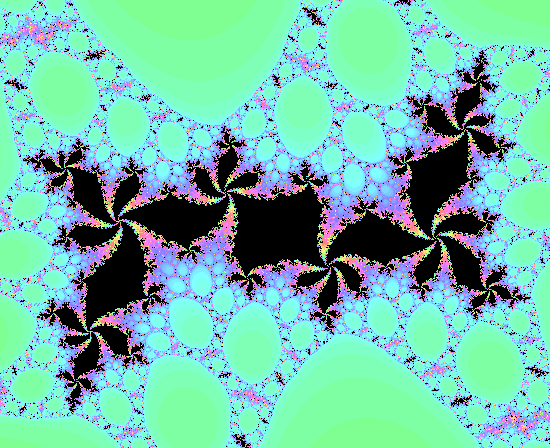}
        \caption[External Angles and Altered Julia on 5-rabbit, $\frac{9}{31}$-bulb, $\FC_-=T1T1$]{Top: A diagram of an altered 9/31-bulb 5-rabbit where $\FC_-=T1T1=\FClabel{41/248}{165/992}{195/992}{49/248}$. Bottom: $\Knot$ where $n=3$, $b=0.03+0.02i$, and $a=0.070054+0.025373i$.
        \label{fig: 9/31 (TR)(TR)}}
    \end{figure}
 \end{example}

\section{Zippers and Zippers with Sidecars}
\label{sec:zippers}

We now turn to the case of baby Julia sets conjugate to quadratic Julia sets with parameter values which lie in \textit{baby} Mandelbrot sets in $\cM$. First, we restrict to the main cardioids of the baby $\cM$'s, and we sometimes split this case in two. More precisely:

\begin{assumptions}
    \label{ass:babyM} Let $\mapR=\Rnab$ satisfy Assumptions~\ref{ass:general}; further, suppose $c$ lies in the main cardioid of a baby Mandelbrot set.   
    
If the baby $\cM$ lies along the real line, call $K(P_c)$ an \textbf{aeroplane}. Else, call it a \textbf{Kokopelli}.
\end{assumptions}

\begin{figure} \centering
            \includegraphics[width=.9\textwidth]{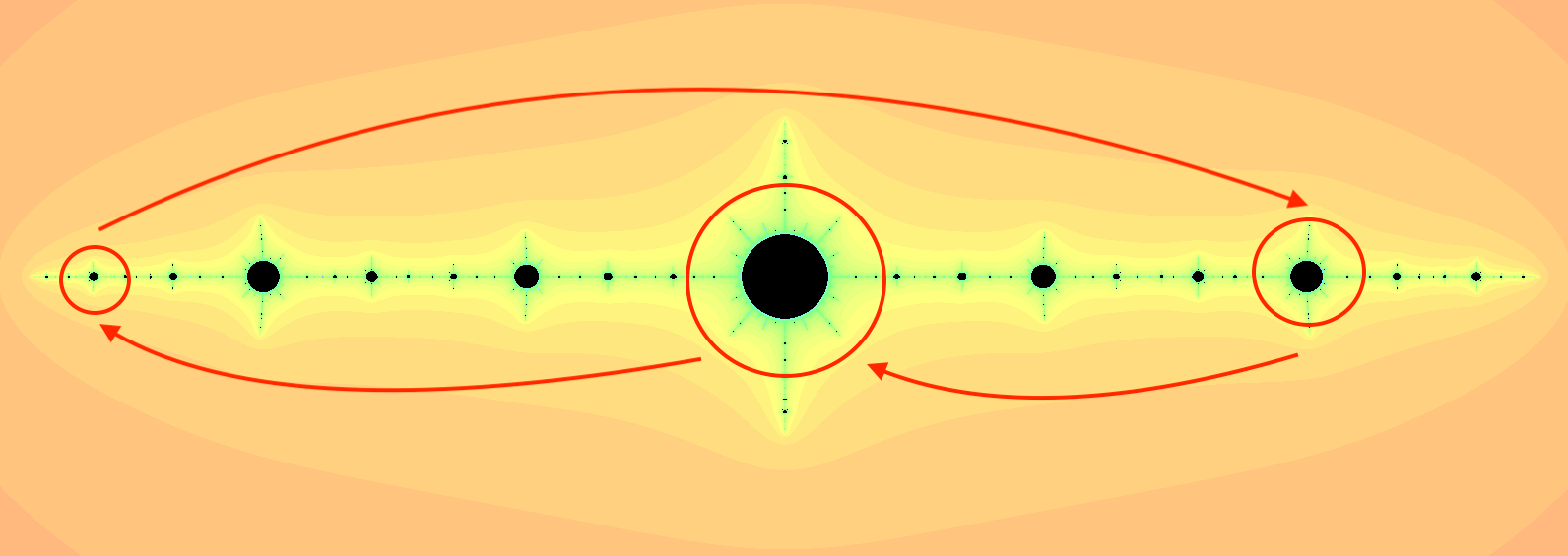}
            \includegraphics[width=.6\textwidth]{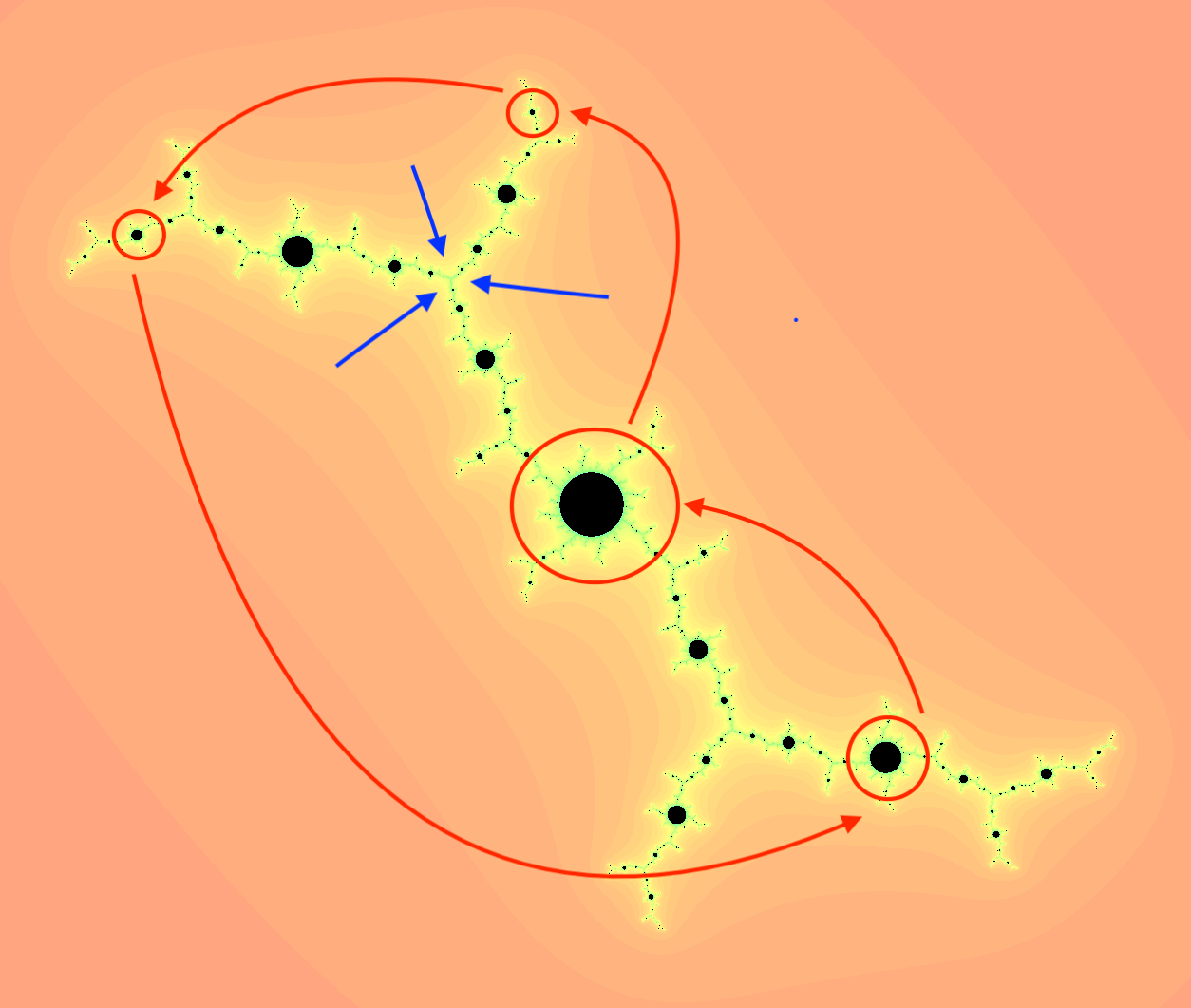}
        \caption[Standard 3-aeroplane and 4-Kokopelli]{ Top: A quadratic 3-aeroplane, $P_c$ for $c\approx -1.755$. Bottom: A quadratic 4-Kokopelli, $P_c$, $c\approx -0.1565 + 1.032i.$; a set of 3 arrows points at a spine junction. In both images, critical orbit Fatou components are circled, and arrows indicate the cycles.
        \label{fig: std aero and koko}}
    \end{figure}

Figure~\ref{fig: std aero and koko} shows a standard 3-aeroplane (attracting orbit period 3) and a 4-Kokopelli (attracting orbit period 4). Reacll that each Fatou component of the filled Julia set of a rabbit is directly adjacent to other Fatou components, so in alterations we described taking ``steps'' through the nearby components. This is not the case for aeroplanes and Kokopellis. Imagine being near the boundary of one of the Fatou components of the filled Julia set of an aeroplane or Kokopelli: when you look out at what Fatou components are nearby, what you see is an infinite chain of Fatou components of size shrinking to zero approaching the boundary of your component. 

\begin{definition} Under Assumptions~\ref{ass:babyM}, 
\begin{enumerate}
    \item[(a)] The \textbf{spine} of an aeroplane or Kokopelli is a tree structure in the filled Julia set which connects any two points of the filled Julia set along a direct 1-dimensional path. The spine runs as ``straight'' as possible through interior Fatou components to get from one entry point, to the center, then out another entry point. 
 
    \item[(b)] A \textbf{spine junction} is a branching point of the spine which does not occur inside of a Fatou component which occur only in Kokopelli Julia sets, not aeroplanes. An example of a spine junction is indicated by three arrows in the Kokopelli of Figure~\ref{fig: std aero and koko}.

\item[(c)] External angles are still defined by a map of $S^1$ onto such a $J(P_c)$, where imagine the path from $0$ to $1$ must always continue along the same side of the spine it was on, working its way around spine junctions. 

Our alteration results relate to locating angles which land on these spines, focusing on angles which land on points at which components meet spines. For both aeroplane and Kokopelli Julia sets, all angles that are identified are identified in pairs (rather than larger sets like 3-rabbits which have angles identified in sets of three, etc.).

When an aeroplane baby Julia set or a Kokopelli baby Julia set is altered, a \textbf{zip} is when all angles within a certain interval which land on a spine swap their identified angles with ones from the landing point symmetric to its landing point.

The name compares the angles splitting from their original identifications to the unzipping of a zipper, where the reidentification is zipping two pieces from separate zippers back together. 

\item[(d)] Along the spine path from $\FC_+$ to $\FC_-$, we consider two components to be \textbf{one zip step away} from each other if the centers of these components can be connected by a piece of spine which passes straight through each component it passes through, entering and leaving such components along the same trajectory, and only making turns at spine junctions. For aeroplanes, this means that single spine pieces visually appear as straight lines, where if two components are connected by spine pieces that meet at a right angle within the component in which they meet, then the components are more than one zip apart. 

This idea carries over for Kokopelli Julia sets, but since Kokopellis have spine junctions that occur outside of components, we allow for any branch to be taken at a spine junction as part of a single zip.

\end{enumerate}
    
\end{definition}

This spine structure makes it significantly more difficult to give the components practical names as we did for $\nu$-rabbits, so we primarily refer to these components by their identified angle names. By convention, we name each component using the identified angles of smallest denominator. For the sake of clarity in this work, we have given components of interest practical names, but suggest that these should not be permanent names outside of the context of the relevant examples.

In aeroplane Julia sets, many of the major spines align with each other in straight lines, which meet each other at right angles inside of components, and this structure allows us to say more about which angles are identified on aeroplanes. Aeroplane Julia sets also have all spines originating from within components. 

For aeroplanes and Kokopellis, since each identification only consists of two angles, we always have that the first angle of each pairing swaps--noting that for angles near $0\equiv1$, the first angle is the one closer to $1$ (non-reduced), and the other is considered to be reduced mod $\mathbb{Z}$ and therefore appears closer to $0$. 

 We can now state the theorem describing aeroplane and Kokopelli alterations. 

\begin{theorem} \label{thm: zip type N}
    Let $F=\Rnab$ satisfy Assumptions~\ref{ass:babyM}.
    
    In all cases, all identifications on $\Jminus$ not mentioned below are unchanged on $\Jnot$. 

    Suppose $\FC_-$ can be reached from $\FC_+$ by traveling through $N$ zip steps along spines of $\Jminus$. Starting with $\FC_+=\FC_0$, let the component at the spine junction where spine $\spire_{i}$ meets spine $\spire_{i+1}$ be labeled as $\FC_i$ with $i=1,\dots,N-1$, letting $\FC_-=\FC_N$. For each $i\in 1,\dots, N$, let the preimage components of $\FC_i$ under $P_c$ be 
    $\FC_i^1$ and $\FC_i^2$
    where the components are listed in the order encountered along $\gamma_-(t)$ as $t$ goes from $0$ to $1$, where $\gamma_- \colon S^1 \to J_-$ is the map defining external angles on $J_-$. As $\FC_i$ is connected to $\FC_{i+1}$ by the spine $\spire_{i+1}$, these preimages will be connected by the preimage of the spine which we will call $\spire_{i+1}^0$, as this is a spine in $\Jnot$. Note there are two sides to $\spire_{i+1}^0$, but they connect at the central component at step $i$ and hence will be considered one piece.
    
    At the $i^{th}$ step, the preimages $\FC_i^1$ and $\FC_i^2$ will lie on opposite sides of the same spine $\spire_i^0$, where the two parts of the spine meet the central component of $\spire_i^0$ at the points with identified angles $c_i^1 \sim d_i^1$ and $c_i^2 \sim d_i^2$, where we assume without loss of generality that $c_i^1<c_i^2$. Let $a_i^2$ be the angle which lands at the point at which $\FC_i^2$ meets $\spire_i^0$. On $\Jminus$, consider any angle $\alpha_i^2 \in [c_i^2, a_i^2]$ where $\gamma_-(\alpha_i^2)$ lies on $\spire_i^0$ and so is identified with another angle $\beta_i^2$. The point at which $\alpha_i^2$ (and $\beta_i^2$) land must be symmetric to a point at which two other angles $\alpha_i^1 \sim \beta_i^1$ are identified. We assume that $\alpha_i^2 < \beta_i^2 < \alpha_i^1 < \beta_i^1$, where $\beta_i^1$ may be reduced mod $\ZZ$. 
    
    On $\Jnot$, all such identifications change to 
    $\alpha_i^1 \sim \beta_i^2$ and $\alpha_i^2 \sim \beta_i^1$;
    that is, given a set of identified angles that land on $\spire_i$, the preimages of those identified angles land on $\spire_i^0$ and swap identified pairs on $\Jnot$.
    This process is repeated for each $i\in 1, \dots, N$.
\end{theorem}

As in Theorem~\ref{thm: rabbit type N}, this is simply an adaptation of Theorem~\ref{thm:typeN}. In this case, angles are again identified in sets of two, but an infinite number of identifications change at once. That being said, as we primarily consider angles that land where components meet spines, we can think of the central component being pinched apart to allow two spines to come together, unzipping the top and bottom parts of each spine from their original identification and reidentifying in opposite halves. Then, as symmetric pairs of components along the spine reach the center, the pairs combine into a new central component, which splits again to allow the next section of spine to come together. If the identifications are arranged in sequence, it is simply the paper fortune teller alteration occurring over and over again in sequence. Thus we provide no formal proof of this theorem.

This Theorem is best understood through the use of examples, which we provide in the following sub-sections; e.g., Example~\ref{ex: 3-aero v_- in 2-7}.

\subsection{Zippers: Aeroplanes}
\label{sec: aero}

Before providing some altered aeorplane examples, we begin by describing the standard aeroplane.

\begin{figure} \centering
    \includegraphics[width=\textwidth]{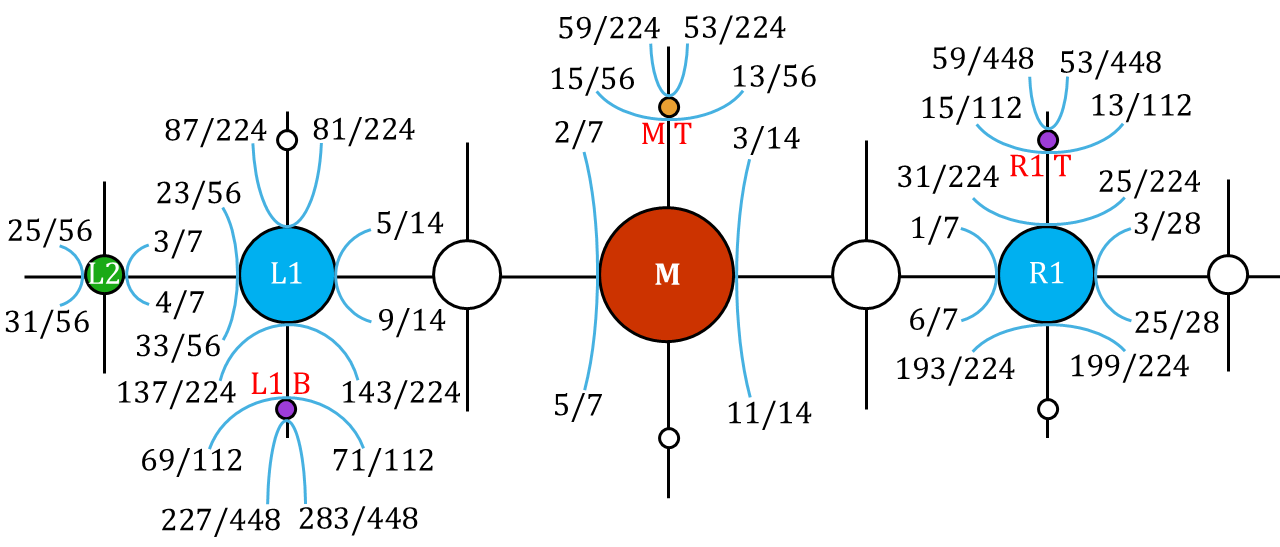}
    \caption[External Angles, standard 3-aeroplane]{A diagram showing the external angles on an unaltered 3-aeroplane\label{fig: 3-aero std}.}
\end{figure}

There are a few situations in which one can easily describe which angles are identified on the aeroplane. We describe them briefly below. Note that while $\nu$-rabbits had $\nu$ components met at each component junction, and hence external angles were identified in groups of $\nu$, to go with their attracting cycle of period $\nu$, here we find angles identified in groups of two, even though all aeroplane Julia sets have attracting cycle with length greater than two. We provide a diagram of the external angles on the 3-aeroplane in Figure~\ref{fig: 3-aero std} to aid in the comprehension of the following results.

We use the following well-known result about external angles in aeroplanes. 

\begin{lemma}
    Consider a quadratic polynomial $P_c$ where $c$ is taken from an aeroplane hyperbolic component of $\mandel$. Recall $\gamma_c \colon S^1 \to J(P_c)$ is the map defining external angles on $J(P_c)$. 
    
    If $\gamma_c(t)=a\in\mathbb{R}$, then $t \sim (1-t)$. That is, if an external angle $t$ lands on a point on the real line, then $t$ is identified with $1-t$. 

    If $\gamma_c(t)$ is a point on a vertical spine of $K(P_c)$, let $t_0 \in [0,t)$ be the greatest angle less than $t$ such that $\gamma_c(t_0)=a_0 \in \mathbb{R}$ and let $t_1 \in (t,1)$ be the least angle greater than $t$ such that $\gamma_c(t_1)=a_1\in\mathbb{R}$. Set $a=(a_0+a_1)/2$. Then $t \sim (a-t)$ mod $\mathbb{Z}$.
\end{lemma}

Values of $c$ for which $K(P_c)$ is a 3-aeroplane come from the main cardioid of the largest baby Mandelbrot set that lies along the real axis in $\mandel$, which can be identified as the landing point of the external ray $\frac{3}{7}$. For such a $\Kplus$, we can use the fact that $M=\FClabel{3/14}{2/7}{5/7}{11/14}$ and $\FC_+=L2=\FClabel{3/7}{25/56}{31/56}{4/7}$. We chose the name $L2$, ``left 2'' for clarity in the examples that follow. These components are colored red and green, respectively, and are part of an attracting 3-cycle with $R1$, colored blue in Figure~\ref{fig: 3-aero std}.

Before examining examples, we provide a corollary to Theorem~\ref{thm: zip type N} specialized for 3-aeroplanes, limiting our scope to where $\FC_-$ is one zip step away from $\FC_+$ along either the real line or the vertical spine running through $\FC_+$, where details about the angle identifications and reidentifications can be made more explicitly. 

\begin{corollary} \label{cor: 3-aero type 1}
    Suppose $\Rnab$ satisfies Assumptions~\ref{ass:babyM} such that $\Jplus$ is a baby 3-aeroplane Julia set.
    Suppose $v_-$ lies in a component $\FC_-$ of $\Jminus$ where $\FC_- \neq \FC_+$. Suppose $\FC_-$ is one zip away from $\FC_+$, which is to say, the centers of $\FC_-$ and $\FC_+$ are connected by a spine $\spire$ which always passes through components at the same trajectory it enters them from, making no turns inside of components.
    Let the two preimage components of $\FC_-=\FC_1$ in $K(P_c)$ be $\FC_1^1$ and $\FC_1^2$, where $\FC_1^1$ is the component encountered first along $\gamma_-(t)$ as $t$ goes from 0 to 1. Recall that $\gamma_-:S^1\to\Jminus$ is the map which assigns external angles to $\Jminus$.

    Suppose $\FC_1^1$ and $\FC_1^2$ both lie along the real axis. Then let $a_2$ be the point at which $\FC_1^1$ meets the spine $\spire_1^0$ which connects the centers of $\FC_1^1$ and $\FC_1^2$. Take $\alpha \in [a_2, \frac{3}{14}]$ to be an angle which lands on the real axis, that is, where $\gamma_-(\alpha) \in \mathbb{R}$, which must then be identified with $1-\alpha$ on $\Jminus$. Then on $\Jminus$, the point at which $\alpha \sim (1-\alpha)$ land must be symmetric to the point at which $(\frac12-\alpha) \sim (\frac12 + \alpha)$ land. On $\Jnot$, these angles will reidentify as $\alpha \sim (\frac12 - \alpha)$ and $(\frac12 + \alpha) \sim (1-\alpha)$. This is true for each $\alpha \in [a_2, \frac{3}{14}]$ which is identified with another angle, including $a_2$ and $\frac{3}{14}$. 
        
    Alternatively, suppose $\FC_1^1$ and $\FC_1^2$ lie along the imaginary axis. Then $\FC_1^1$ meets the spine $\spire_1^0$ which connects their centers at a point $a_1$. Take $\alpha \in [\frac{25}{112}, a_1]$ to be an angle which lands on the imaginary axis, that is, where $\gamma_-(\alpha) \in \mathbb{R}i$. Then on $\Jminus$, we have $\alpha \sim (\frac12-\alpha)$ and $(\frac12+\alpha) \sim (1-\alpha)$. On $\Jnot$, we will instead have $\alpha \sim (1-\alpha)$ and $(\frac12-\alpha) \sim (\frac12 + \alpha)$. This is true for each such $\alpha$. 
\end{corollary}
The proof is a straightforward application of Theorem~\ref{thm: zip type N} to 3-aeroplanes, whose unaltered external angles are documented in Figure~\ref{fig: 3-aero std}.

\begin{figure} \centering
        \includegraphics[width=.4\textwidth]{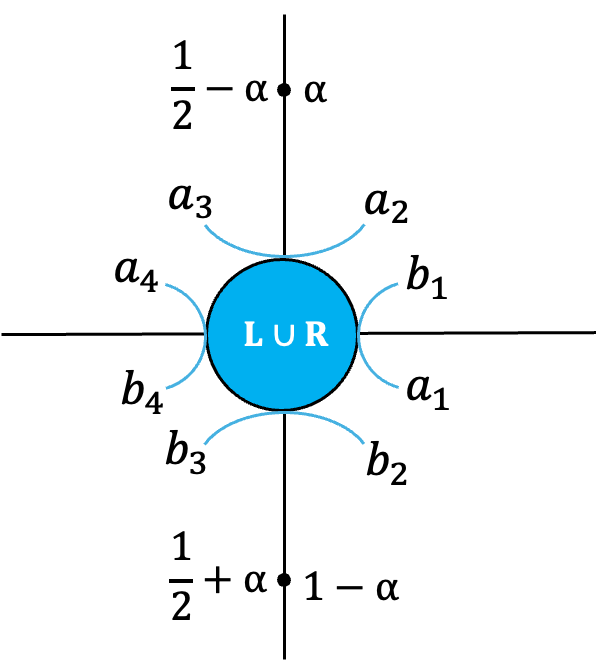}
    \includegraphics[width=.95\textwidth]{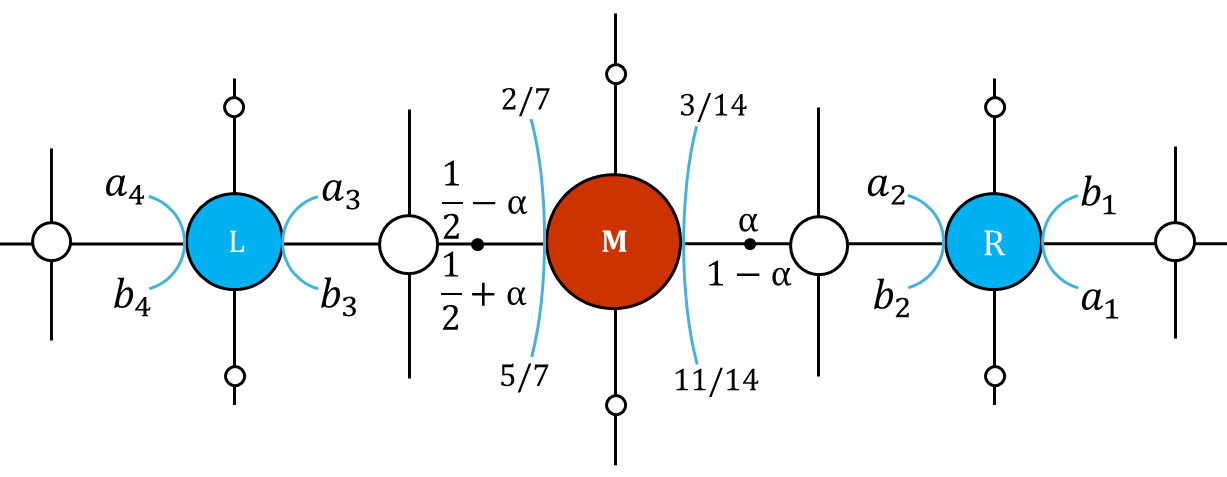}
    \caption[External Angle Identification changes for $\FC_-$ on $\RR$]{A diagram showing the angles in the 3-aeroplane which split and reidentify when the preimages of $\FC_-$ lie on the real line. Top: $\Jnot$, Bottom: $\Jminus$. \label{fig: 3-aero cor diagram}}
\end{figure}

As a way to envision the angle identification changes listed in the first case, 
imagine that the left and right halves of the spine $\spire$ running along the real axis in a 3-aeroplane are each a zipper, where
the top half is zipped to the bottom half for both the left and right sides. On $\Jnot$, we unzip and rezip so that the two top sides are zipped to each other, and the same on the bottom. However, we are not fully unzipping the left and right zippers, but only the portion that lies between the two preimages of $\FC_-$ in $K(P_c)$, leaving the remainder of the left and right halves unchanged. See Figure~\ref{fig: 3-aero cor diagram} for a diagram of this alteration, where the components labeled $R$ and $L$ represent $\FC_1^1$ and $\FC_1^2$, respectively.

Now we provide a pair of examples to demonstrate Theorem~\ref{thm: zip type N} and Corollary~\ref{cor: 3-aero type 1}: a type $N=1$ example, then a type $N=2$ which builds upon it. 

\begin{example} \label{ex: 3-aero v_- in 2-7}
    First, consider the case where $v_-$ lies in $\FC_-=M=\FClabel{3/14}{2/7}{4/7}{11/14}$. This component lies on the spine $\spire$ which runs along the real axis, so it can be reached from $\FC_+$ in one step. The preimages of $\FC_-=\FC_1$ in $K(P_c)$ are $\FC_1^1=R1=\FClabel{3/28}{1/7}{6/7}{25/28}$ and $\FC_1^2=L1=\FClabel{5/14}{11/28}{17/28}{9/14}$, which are colored blue in Figure~\ref{fig: 3-aero std}. To combine these two components into one, we need to unzip everything between them and $M$, which coincidentally in this case happens to be $\FC_-$, since it lies in the center of the spine $\spire^0$ which connects $\FC_1^1$ and $\FC_1^2$. Therefore we consider each angle $\alpha \in [1/7, 3/14]$ which lands on $\spire^0$. In this case, since $\spire^0$ terminates in 0 and $\frac12$, we have on $\Jminus$ that $\alpha \sim (1-\alpha)$ and $(\frac12 - \alpha) \sim (\frac12 + \alpha)$. On $\Jnot$, these pairings change to $(\frac12+\alpha) \sim (1-\alpha)$ and $\alpha \sim (\frac12 - \alpha)$ for each $\alpha$. This keeps all spines which decorate components along $\spire^0$ intact, while taking each set of symmetric components between $\FC_1^1$ and $\FC_1^2$, unzipping the top and bottom halves of each one from each other, and zipping back together the two tops and the two bottoms, until finally we split apart $\frac{5}{14} \sim \frac{9}{14}$ and $\frac67 \sim \frac17$ and reidentify these angles as $\frac67 \sim \frac{9}{14}$ and $\frac{5}{14} \sim \frac{1}{7}$. In the diagram included in Figure~\ref{fig: 3-aero v_- in 2-7}, one can see these newly identified angles on the bottom and top of the new central component, respectively.  \qed

    \begin{figure} \centering
            \includegraphics[width=\textwidth]{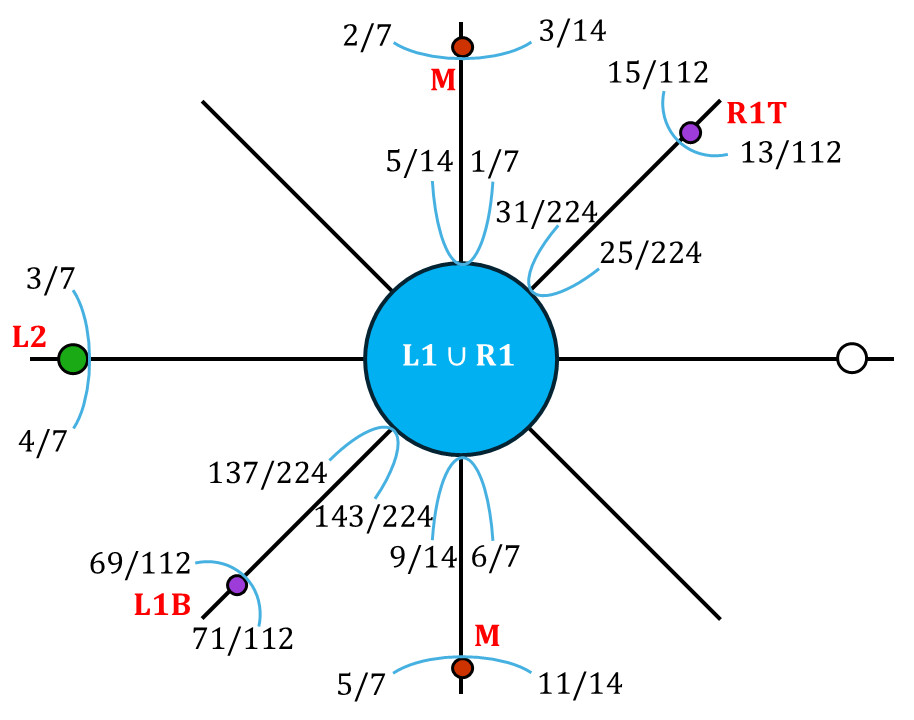}
            \includegraphics[width=.45\textwidth]{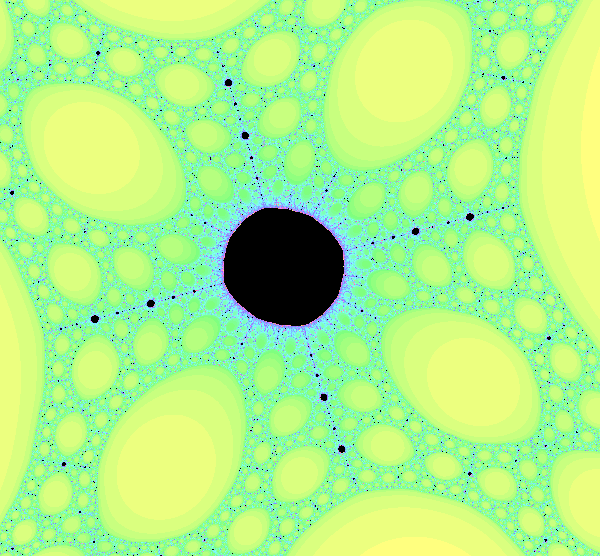}
        \caption[External Angles and Altered Julia set, 3-aeroplane, $\FC_-=M$ ]
        {Top: A diagram of an altered 3-aeroplane where $\FC_-=\FClabel{3/14}{2/7}{5/7}{11/14}$. Bottom: $\Knot$ where $n=3$, $b=0.02+0.03i$, and $a=0.01949-0.01126i$.\label{fig: 3-aero v_- in 2-7}}
    \end{figure}
\end{example}

Next is a Type $N=2$ example which builds off of the prior example.

\begin{example} \label{ex: 3-aero v_- in 13-56}
    Suppose that $v_-$ lies in $\FC_-=MT=\FClabel{13/56}{53/224}{59/224}{15/56}$. This component lies above $M$ in $K(P_c)$, so we see that it is reached from $\FC_+=\FC_0$ by first traveling along the real axis to $\FC_1=M=\FClabel{3/14}{2/7}{5/7}{11/14}$, and then traveling along the imaginary axis to $\FC_-=\FC_2$, thus being called $MT$ for ``main top''. Hence this case will be considered to be of Type $N=2$. Observe that the two preimages of $\FC_-$ under $P_c$ are $\FC_2^1=R1T=\FClabel{13/112}{53/448}{59/448}{15/112}$ and $\FC_2^2=L1B=\FClabel{69/112}{277/448}{283/448}{71/112}$, colored purple, which do not lie on the same spine in the unaltered 3-aeroplane as shown in Figure~\ref{fig: 3-aero std}. However, if we consider the changes made to angle identifications in Example~\ref{ex: 3-aero v_- in 2-7} as a first or intermediate step, then the second step can be considered further alterations to this already altered 3-aeroplane, and in Figure~\ref{fig: 3-aero v_- in 2-7}, we see that $\FC_2^1$ and $\FC_2^2$ do share a spine $\spire_2^0$ in this intermediate step. $\spire_2^0$ terminates in $\frac{7}{56}$ and $\frac{35}{56}$ and meets its central component at the points $\frac{25}{224} \sim \frac{31}{224}$ and $\frac{137}{224} \sim \frac{143}{224}$, so we see it is any angle $\alpha$ which lands on $\spire_2^0$ within the interval $[\frac{25}{224}, \frac{13}{112}]$ that need to be split and reidentified. On this intermediate step and on $\Jminus$, these angles present in the identifications $\alpha \sim (\frac14 - \alpha)$ and $(\frac12 + \alpha) \sim (\frac34-\alpha)$. On $\Jnot$, these identifications change to $(\frac12 + \alpha) \sim (\frac14 - \alpha)$ and $\alpha \sim (\frac34-\alpha)$. It is both sets of changes that are present on $\Jnot$. A diagram of the new identified angles is given in Figure~\ref{fig: 3-aero v_- in 13-56}. \qed

    \begin{figure} \centering
            \includegraphics[width=0.8\textwidth]{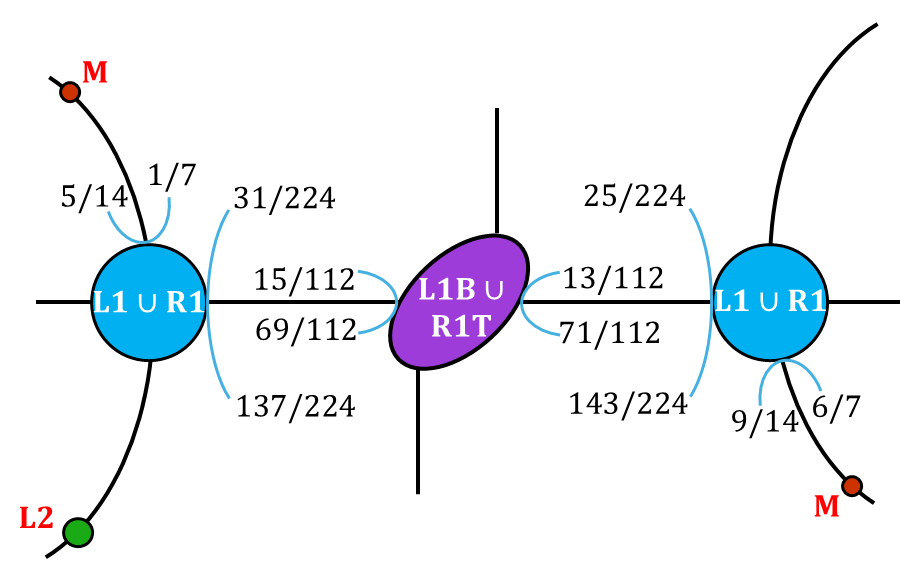}
            \includegraphics[width=.5\textwidth]{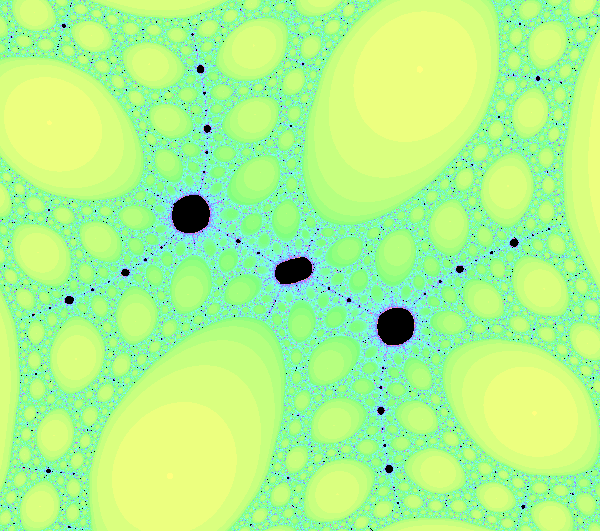}
        \caption[External Angles and Altered $J$, 3-aeroplane, 2 steps: $M$ then up]{Top: A diagram of an altered 3-aeroplane where $\FC_-=MT=\FClabel{13/56}{53/224}{59/224}{15/56}$. Bottom: $\Knot$ where $n=3$, $b=0.02+0.03i$, and $a=0.019491-0.011214i$. \label{fig: 3-aero v_- in 13-56}}
    \end{figure}
 \end{example}

\begin{figure}
    \includegraphics[width=\textwidth]{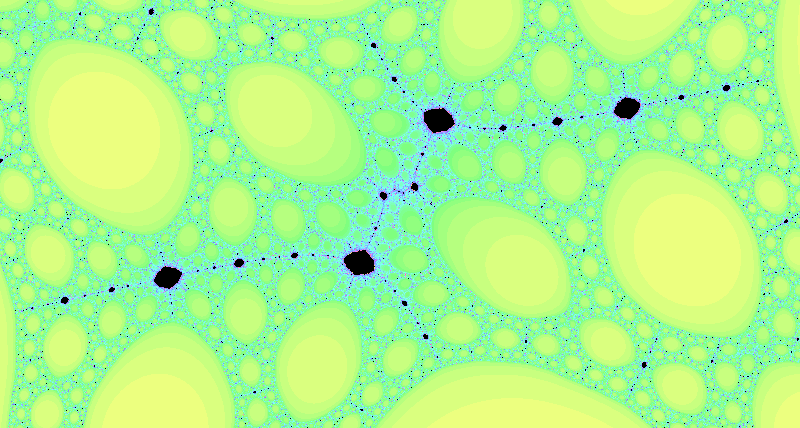}
    \caption{An altered 3-aeroplane of Type $N=4$.\label{fig: 3-aero N4}}
\end{figure}

We conclude this section with a visual of an altered 3-aeroplane that appears to be of Type $N=4$ in Figure~\ref{fig: 3-aero N4}. 
Aeroplanes with longer periodic cycles do not behave differently, so we do not include any examples here, but an interested reader can see \cite{Brouwer-Thesis} for some 4-aeroplane examples.


\subsection{Zippers with Spine Sidecars: Kokopellis}
Next, we proceed to analyze alterations of  Kokopelli Julia sets. The quadratic Julia sets typically referred to as Kokopellis are spawned from the main cardioid of the largest baby Mandelbrot set connected by a spine to the upper 3-rabbit ($\frac{1}{7}$-) bulb, which gives them an attracting cycle of period 4, but we expand our category to include Julia sets spawned from any main cardioid of a baby $\mandel$ that does not lie along the real line. Like aeroplanes, Kokopellis are comprised of components connected by spines, but in this case, the spines do not fall primarily in straight lines, and the attracting cycle has a less obvious pattern as the components involved are positioned off of the primary spines. Furthermore, the spines occasionally branch off, so that three pieces of spine are present at each spine junction.

Because of the spine junctions present in the Kokopelli Julia set, we find that the alterations of these baby Julia sets follow a similar zipper style of splitting and reidentifying as occurred in the aeroplane Julia sets, but with the addition of the ``sidecars'' that were present in $\nu$-rabbits, where in this case the sidecars are an entire section of spine that does not split and reidentify.

\begin{definition} \label{defn:general sidecar}
For Kokopelli Julia sets, we expand the use of the term \textbf{sidecar} (beyond Definition~\ref{defn:sidecar}) to refer to a section of Julia set which ends on a spine junction involved in a zip, but the angles of this section are not altered in the zip. 
\end{definition}

We describe this phenomena more precisely in the following corollary of Theorem~\ref{thm: zip type N}, specific to Kokopelli Julia sets of Type $N=1$.

\begin{corollary}\label{cor: koko N=1} (Kokopelli, Type $N=1$) 
    Suppose $\Rnab$ satisfies Assumptions~\ref{ass:babyM}; further, assume $\Jplus$ is a baby Kokopelli Julia set. Suppose $v_-$ lies in a component $\FC_-\neq\FC_+$ of a preimage of $\Kplus$, where the component of $K(P_c)$ that is identified by the same angles as $\FC_-$ is reachable from $\FC_+$ by traveling along a single spine. Recall that $\gamma_-:S^1\to\Jminus$ is the map which assigns external angles to $\Jminus$.

    Let $\FC_1^1$ and $\FC_1^2$
    be the preimages of $\FC_-$ under $P_c$, where $\FC_1^1$ is the first of these components encountered along $\gamma_-(t)$ as $t$ goes from 0 to 1. Let $\spire_1^0$ be the spine that connects $\FC_1^1$ and $\FC_1^2$. Let $a_1$ be the point at which $\FC_1^1$ meets $\spire_1^0$.

    Note that only spine pieces that contribute to connecting $\FC_1^1$ and $\FC_2^1$ will count as part of $\spire_1^0$, where as sidecar spines that connect to $\spire_1^0$ at spine junctions or spines that originate from components that lie along $\spire_1^0$ will not. 

    Let $\alpha_1$ be an external angle that lands on $\spire_1^0$, which is identified with another angle such that $\alpha_1 \sim \beta_1$. Let $\gamma_-(\alpha_2)$ be rotationally symmetric to $\gamma_-(\alpha_1)$, and let $\alpha_2 \sim \beta_2$ where $\alpha_2$ is chosen so that either $\beta_1$ or $\beta_2$ has value between $\alpha_1$ and $\alpha_2$ when arranged in increasing order. Then on $\Jnot$, these identifications are changed to $\alpha_1 \sim \beta_2$ and $\alpha_2 \sim \beta_1$. Any angle that does not land on $\spire_1^0$ in $\Jminus$ has unaltered identifications on $\Jnot$.
\end{corollary}

\textbf{Counting zip steps.} By comparison to the case stated above, the aeroplane case seems simpler and can be stated in more detail because the structure of the aeroplane seems to suggest that a single spine is any spine pieces that lay in a straight line. Thus, we included the simplified statement of Corollary~\ref{cor: 3-aero type 1} for ease of reading and understanding. However, alterations to baby aeroplane Julia sets are combinatorially no different than the case presented by Kokopelli Julia sets, in which many spines meet at spine junctions. One simply needs to carefully identify which pieces of spine connect the preimages which need to be combined. We find that any subsequent spine may be chosen at a spine junction while still constituting a single zip step. To be more than one step away, the spine path between the preimage components must turn inside of a component rather than continuing out of the component along a spine in the same trajectory at which the component was entered. In this way, the ``in the same line'' idea of what constitutes a spine is carried over from the aeroplane example, although several spine pieces lined up end to end still count as one. State simply, turns that occur at a spine junction count as part of one step, and turns that occur inside of a component count as distinct steps.

\medskip

\begin{figure} \centering
    \includegraphics[width=\textwidth]{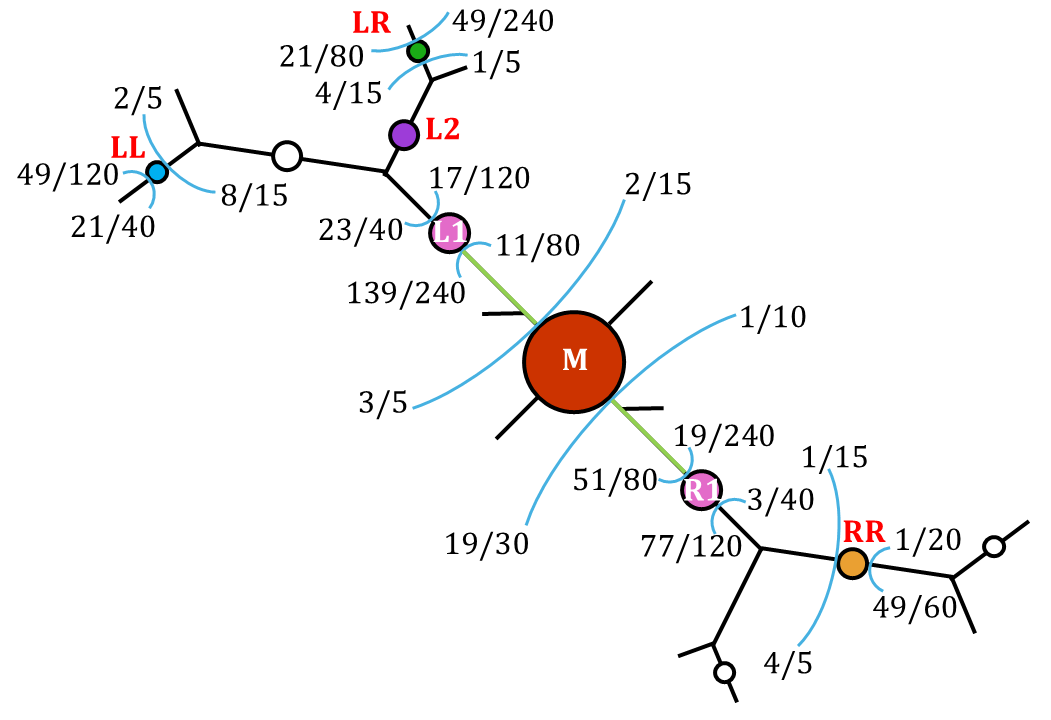}
    \caption{A diagram showing the external angles on an unaltered Kokopelli Julia set.\label{fig: koko std}}
\end{figure}

In the following, it will be useful to note that the critical value component of the standard 4-cycle Kokopelli is 
\\ $\FC_+=LR=\FClabel{1/5}{49/240}{21/80}{4/15}$, and so the period 4 cycle of $\frac15 \to \frac25 \to \frac45 \to \frac35$ is the key focus of this Julia set. We also use the fact that here, $M=\FClabel{1/10}{2/15}{3/5}{19/30}$. A diagram detailing the external angles on an unaltered Kokopelli Julia set is given in Figure~\ref{fig: koko std}. Again, the component names given here are simply for ease of reading in the examples that follow.

The examples we worked through for this case are not as varied as for some previous sections. This is in part due to the difficulty of determining what constituted a Type $N=2$ example, and after the determination was reached, to the highly precise nature of the components involved. Only one Type $N=2$ example was included in \cite{Brouwer-Thesis} because any additional examples of this type dealt with external angles with denominators of at least five figures. Below, we simply show one Type $N=1$ example. The main point of interest of this section is to note that the spine junctions aren't an issue and become the sidecars in alterations, but otherwise, altered Kokopelli Julia sets behave the same as altered aeroplanes.

\begin{example} \label{ex: koko v_- in 3-20}
    As a simplest example, let $\FC_-=L2=\FClabel{3/20}{19/120}{11/40}{17/60}$. The two preimages of $\FC_-=\FC_1$ under $P_c$ are $\FC_1^1=R1=\FClabel{3/40}{19/240}{51/80}{77/120}$ and $\FC_1^2=L1=\FClabel{11/80}{17/120}{23/40}{139/240}$, colored pink in Figure~\ref{fig: koko std}, which appear on the unaltered Kokopelli Julia set as the largest components on either side of the central component lying along the central spine piece. It is this primary spine $\spire_1^0$ that needs to unzip and reidentify, which is why the unaltered $\spire_1^0$ was colored light green in Figure~\ref{fig: koko std}. However, it is not as easy as in the aeroplane case to describe the exact angles that are identified on $\Jminus$ as opposed to $\Jnot$, which is why the description in Corollary~\ref{cor: koko N=1} is written as it was. Taken in symmetric pairs, the starting angles that land on $\spire_1^0$ are a subset of the interval $[\frac{19}{240}, \frac{1}{10}]$, ignoring the spine sidecars such as the two left in black along the light green spine in Figure~\ref{fig: koko std}.
    The first step of the splitting involves pinching $M$ apart into $\FClabel{1/10}{49/480}{21/160}{2/15}$ and $\FClabel{3/5}{289/480}{101/160}{19/30}$, then unzipping the two sides of the existing spines and rezipping them to each other. When the spine sidecars mentioned previously are reached, they stay intact, moving with their side of the spine, as can be seen on the altered figure where the black spines share spine junctions with light green spines. On $\Jnot$, we see the altered spines extending out to the sides of the new central component $\FClabel{3/40}{17/120}{23/40}{77/120}$. In fact, one can see that everything on the side spines protruding from this new central component are where angle identifications have been altered, except for the spine sidecars in black which are unaltered. A diagram detailing these changes is given in Figure~\ref{fig: koko v_- in 3-20}, alongside a computer generated image of the altered baby Julia set. 
    \qed

    \begin{figure} \centering
              \includegraphics[width=.65\textwidth]{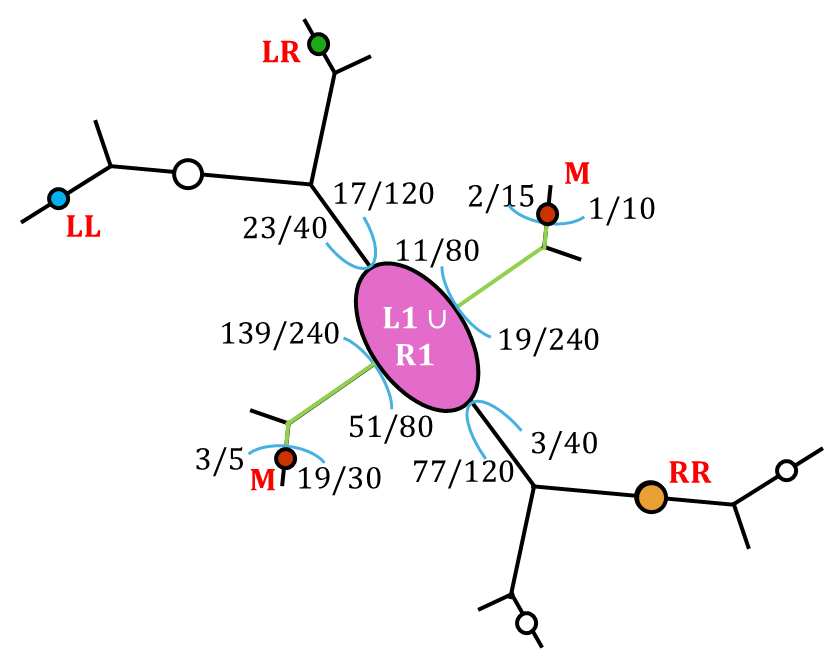}
            \includegraphics[width=.65\textwidth]{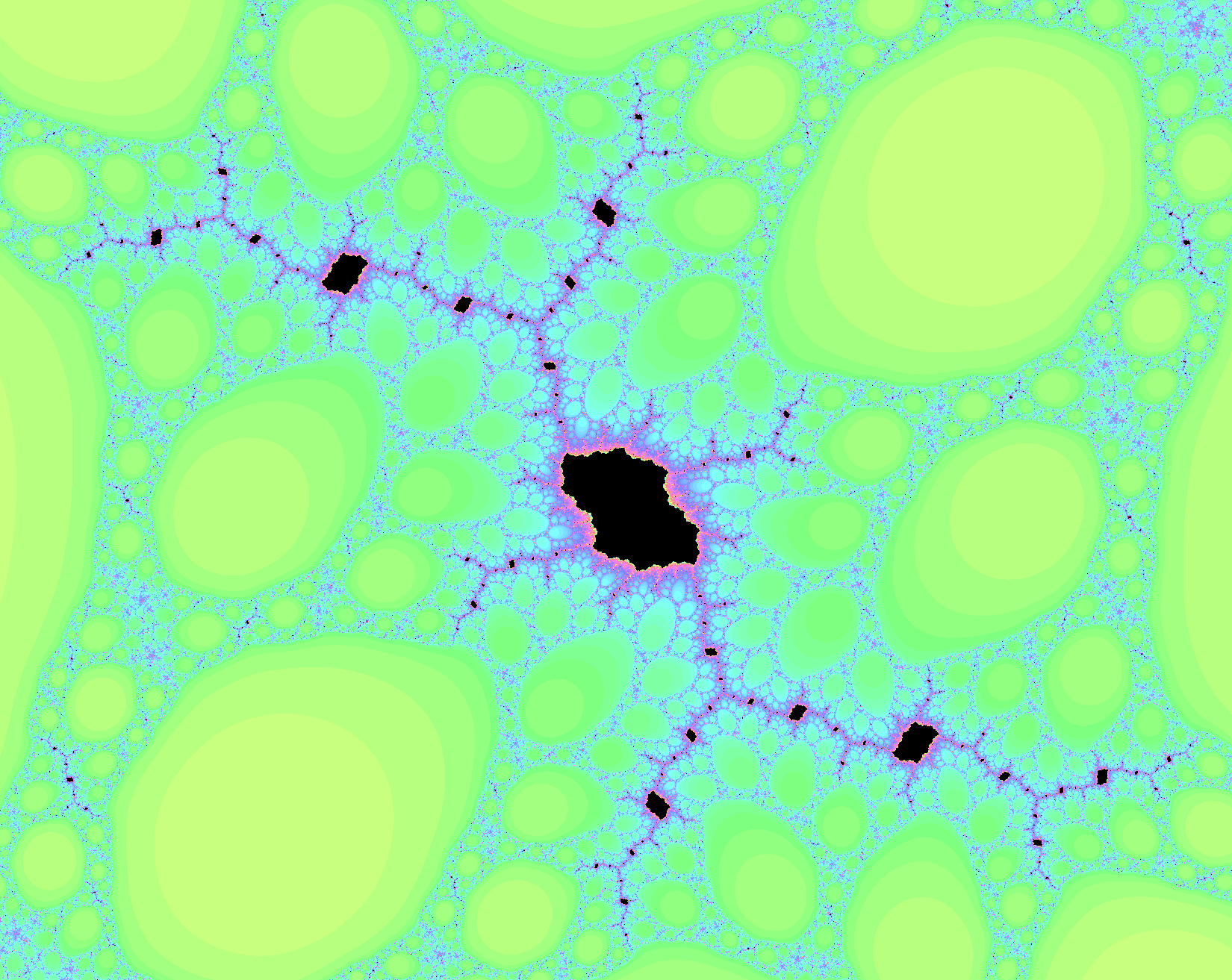}
        \caption[External Angles and Altered $J$ for 4-Kokopelli of Example~\ref{ex: koko v_- in 3-20}]{Top: A diagram showing the identified angles on an altered Kokopelli baby Julia set where $\FC_-=\FClabel{3/20}{19/120}{11/40}{17/60}$. Bottom: $\Knot$ where $n=3$, $b=0.02+0.02i$, and $a=0.084593+.061765i$. \label{fig: koko v_- in 3-20}}
    \end{figure}
 \end{example}

\begin{figure}
    \includegraphics[width=.8\textwidth]{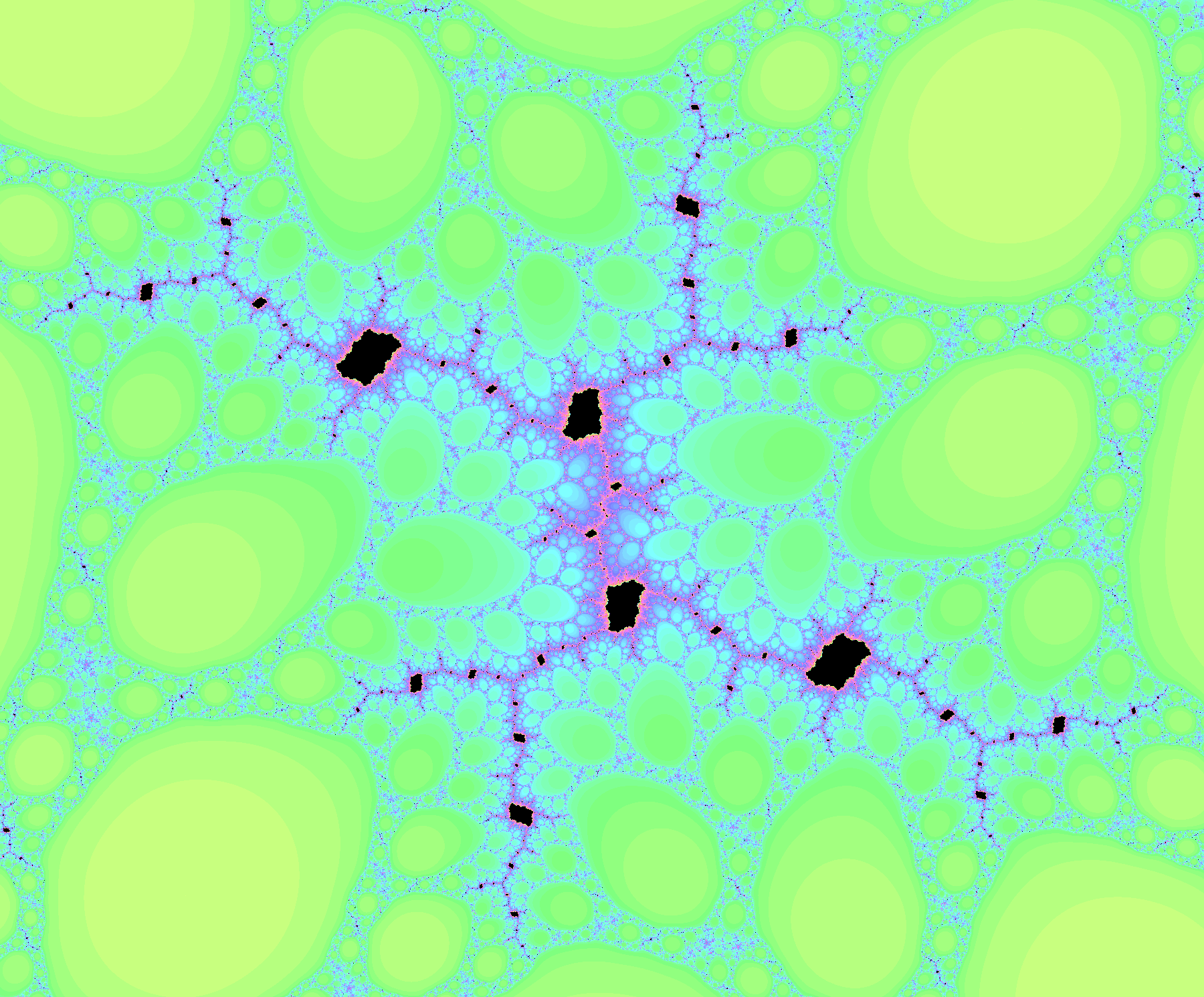}
    \caption{An altered Kokopelli baby Julia set of Type $N=3$\label{fig: koko N=3}.}
\end{figure}

\begin{figure}
    \includegraphics[width=.6\textwidth]{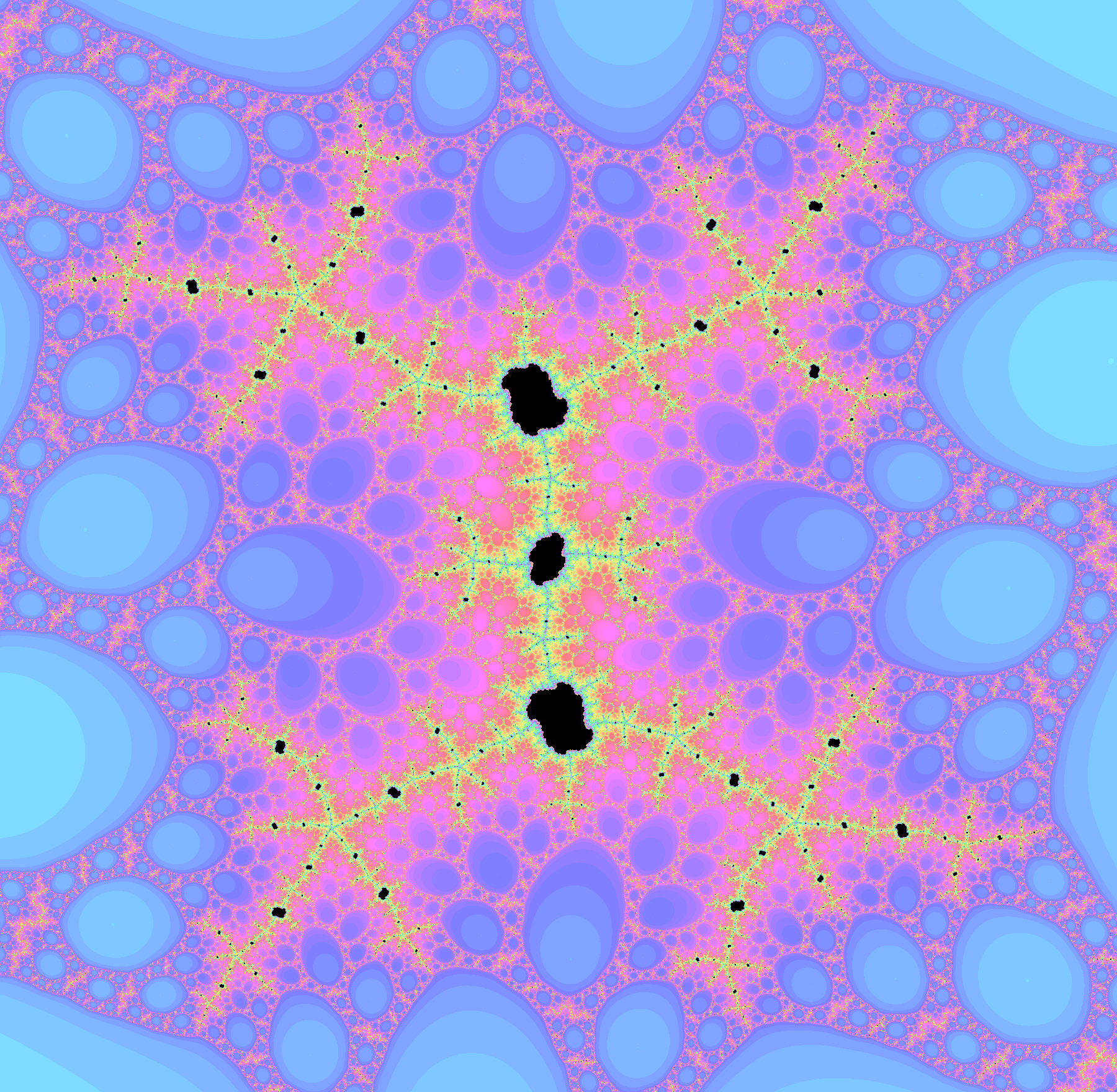}
    \caption{$\Knot$ where $n=3$, $b=0.02+0.02i$, $a=0.06740764+0.028731855i$.\label{fig: 6-koko Type N=2}}
\end{figure}

We close this section with a couple of pictures of more complicated cases, but do not provide details on the angle identification changes.
Figure~\ref{fig: koko N=3} appears to be a Type $N=3$ altered 4-Kokopelli baby Julia set,
and Figure~\ref{fig: 6-koko Type N=2} appears to be an altered baby 6-Kokopelli Julia set of Type $N=2$.
Compare its shape with that of similar examples in Figures~\ref{fig: 3-rabbit upper v_- in M(RT) altered} and~\ref{fig: 3-aero v_- in 13-56}.

\section{Future Work}
\label{sec:future}

To finish the catalog of potential Julia alterations, our next case study involveS alterations of the external angles on all remaining hyperbolic components of $\mandel$. The components we have not yet discussed are any non-primary decorations of the main cardioid of $\mandel$, or any decorations of main cardioids of baby Mandelbrot sets. Each of these components spawns a different kind of quadratic Julia set with different external angle identifications, but all of them share the quality in that they take the general form of a Julia set whose alterations we have already discussed, with each component replaced with an entire Julia set. 

Other future work could be to interpret the angle changes in the language of Thurston's lamination diagrams, or examine the parabolic case.


\bibliographystyle{alpha} 

\end{document}